\documentclass[11pt]{article}

\usepackage{tikz}
\usepackage{subfigure}
\usepackage{hyperref}
\usepackage[text={450pt,600pt}]{geometry}
\usepackage{titlesec}
\usepackage{amssymb}
\usepackage{xcolor}
\usepackage{amsmath,amssymb,amscd,mathrsfs}
\usepackage{color}
\usepackage{amssymb}
\usepackage{amsmath}
\usepackage{amscd}
\usepackage{ntheorem}
\usepackage{xcolor}

\DeclareMathOperator{\Prog}{Prog}

 \newcommand{\sgn}{\operatorname{sgn}}
 
\usepackage{cleveref}

\hypersetup{
  colorlinks   = true,    
  urlcolor     = blue,    
  linkcolor    =  blue,  
  citecolor    = blue      
}
\newcommand*{\qed}{\par\vspace{5mm} \hfill\ensuremath{\Box}}
\usepackage[latin1]{inputenc}   
\usepackage[T1]{fontenc}   
\usepackage{mathtools}
\DeclarePairedDelimiterX{\norm}[1]{\lVert}{\rVert}{#1}

\begin{document}

\renewcommand\thesubsection{\arabic {subsection}}
\newtheorem{thm}{Theorem} [section]
  \newtheorem{Conjecture}[thm]{Conjecture}
   \newtheorem{Lemma}[thm]{Lemma} 
\newtheorem{prop}[thm]{Proposition}
 \newtheorem{Corollary}[thm]{Corollary}
     \newtheorem{ex}{Example}
    \theorembodyfont{\normalfont}
  \newtheorem{rem}{Remark} 
     \newtheorem{as}{Assumption}[section]
\renewcommand{\theas}{{{(A.\arabic{as}})}}
\newtheorem{defi}{Definition}[section]
  \renewcommand{\thedefi}{\arabic{defi}.}
  \renewcommand{\therem}{\arabic{rem}.}
   \theorembodyfont{\normalfont}\theoremstyle{nonumberplain}  
     \newtheorem{Proof}{Proof.}
  \newtheorem{remm}{Remark.} 

  \renewcommand{\thethm}{\arabic{thm}}%

\newenvironment{thmbis}[1]
  {\renewcommand{\thethm}{\ref{#1}$'$}%
   \addtocounter{thm}{-1}%
   \begin{thm}}
  {\end{thm}}

\newenvironment{asp}[1]
 {\addtocounter{as}{-1}%
  \renewcommand{\theas}{(A.\arabic{as}$'$)}%
   \begin{as}}
  {\end{as}}

\newenvironment{aspp}[1]
 {\addtocounter{as}{-1}%
  \renewcommand{\theas}{(A.\arabic{as}$''$)}%
   \begin{as}}
  {\end{as}}

\makeatletter
\newcommand{\neutralize}[1]{\expandafter\let\csname c@#1\endcsname\count@}
\makeatother
\title{\bf $\mathbb{L}^p$ Solutions of Quadratic BSDEs}
\author{Hanlin Yang\thanks{{Department of Banking and Finance, Universit\"at Z\"
urich, Plattenstrasse 22, CH-8006 Z\"urich, Switzerland, e-mail: hanlin.yang@uzh.ch. This paper is part of the author's master thesis.}}} 
\date{\today}
\maketitle
\abstract{
We  study a general class of quadratic BSDEs with   terminal value in   $\mathbb{L}^p$ for $p>1$. First of all, we give an $\mathbb{L}^p$-type estimate and existence result. Under the additional assumption of  monotonicity and convexity, we derive the comparison theorem, uniqueness and stability result via $\theta$-technique (Briand and Hu \cite{BH2008}). The assumptions employed throughout this paper are rather weak and extend the quadratic BSDE literature. 
Finally,  a probabilistic representation for the viscosity solution to the associated quadratic PDEs is given. }
\subsection{Introduction}
In this paper, we are concerned with $\mathbb{R}$-valued backward stochastic differential equations (BSDEs)
\begin{align}
Y_t = \xi + \int_t^T F(s, Y_s, Z_s)ds - \int_t^T Z_s dW_s,  \label{bsde}
\end{align}
where 
 the generator $F$ is continuous and satisfies $\mathbb{P}$-a.s. for all $(t, y, z) \in [0, T]\times\mathbb{R}\times\mathbb{R}^d$, 
\begin{align}
\sgn(y)F(t, y, z)&\leq
\alpha_t + \beta|y|+ \gamma |z|+ f(|y|)|z|^2, \nonumber \\
|F(t, y, z)|&\leq \alpha_t + \varphi(|y|) + \gamma |z| + f(|y|)|z|^2, 
\label{lq0}
\end{align}
for an $\mathbb{R}^+$-valued progressively measurable process $\alpha$, $\beta\in\mathbb{R}, \gamma \geq 0, $   a function
$f(|\cdot|): \mathbb{R}\rightarrow \mathbb{R}^+$  which is integrable and bounded on any compact subset of $\mathbb{R}$, and a continuous nondecreasing function $\varphi:\mathbb{R}^+\rightarrow\mathbb{R}^+$. 
A solution to (\ref{bsde}) is a process $(Y, Z)$ adapted to the filtration generated by the Brownian motion $W$ such that (\ref{bsde}) holds $\mathbb{P}$-a.s. for all $t\in[0, T]$.
We emphasize that, unlike the  quadratic BSDEs studied by Briand and Hu \cite{BH2006}, \cite{BH2008}, 
the quadratic growth in our study takes the  form  $f(|y|)|z|^2$. 
 Moreover, we assume that the terminal value $\xi$ and $\int_0^T \alpha_s ds$ belong to $\mathbb{L}^p$ for a certain $p>1$.

Let us recall that, 
quadratic BSDEs
 are first studied by Kobylanski \cite{K2000}, where 
 existence, uniqueness, comparison theorem and  monotone stability  for bounded solutions are obtained. 
 Proving the existence of a solution consists in constructing a monotone sequence of bounded solutions of better-known BSDEs and then passing the limit. The underlying machinery of this procedure is called the monotone stability of quadratic BSDEs. 
  Later, Briand and Hu \cite{BH2006}, \cite{BH2008}
extend the existence result  by assuming that the terminal value has  exponential moments integrability. Recently, Bahlali et al \cite{BEO2014} constructs a solution to quadratic BSDEs with its terminal value in $\mathbb{L}^2$ and a  generator 
  satisfying
\begin{align*}
|F(t, y, z)|\leq \alpha + \beta|y| + \gamma |z| + f(|y|)|z|^2, 
\end{align*}
for some $\alpha, \beta, \gamma \geq 0$.
However, as to the uniqueness of a solution, only purely quadratic BSDEs are studied. 

There are two lines of studies on the uniqueness of a solution to quadratic BSDEs. When the terminal value is bounded, one crucial feature is that $\int_0^\cdot Z_sdW_s$ is  a BMO martingale. This observation, combined with a local Lipschitz condition, can be used to prove a uniqueness result; see, e.g., \cite{HIM2005}, \cite{MS2005}, \cite{M2009}, \cite{BE2013}. However,
 $\int_0^\cdot Z_sdW_s$ is in general not a BMO martingale
if the terminal value is unbounded.  
  Nevertheless one can also obtain a uniqueness result, by relying on a 
  convexity condition which proves to be convenient to treat  the quadratic generators; see \cite{BH2008}, \cite{HS2011}, \cite{DHR2011}, etc.

The first contribution of this paper is to study  an existence result given  (\ref{lq0}) and a terminal value in $\mathbb{L}^p$ for a certain $p>1$. We first briefly present the motivations  to assume (\ref{lq0}). 
Among the literature on non-quadratic BSDEs, assumptions  of this type are quite convenient  to obtain the a priori estimates; see, e.g., \cite{BC2000}, \cite{B2003}, \cite{BLS2007}.
It turns out that the existence and monotone stability of bounded solutions can also be adapted to quadratic BSDEs with a growth  of this type.  This is stated in 
 Briand and Hu \cite{BH2008}, which assumes that
\begin{align*}
\sgn(y)F(t, y, z) &\leq \alpha_t + \beta|y| + \eta|z|^2, \\
|F(t, y, z)|&\leq \alpha_t + \varphi(|y|) + \eta|z|^2.
\end{align*}
  The proof is merely a slight modification of  Kobylanski \cite{K2000}. In parallel with these works, we prove an existence result under (\ref{lq0}).  In the first step, we derive a $\mathbb{L}^p$-type estimate for quadratic BSDEs, by adapting the method developed by Briand et al \cite{B2003}. To construct a solution, we use a combination of  the localization procedure   developed by Briand  and Hu \cite{BH2006} and the monotone stability result.

Another contribution is to address the question of uniqueness. 
In the spirit of Briand and Hu \cite{BH2008}, we prove comparison theorem, uniqueness and a stability result via 
$\theta$-technique
 under a monotonicity and convexity assumption. 
It turns out that, our results of existence and uniqueness, not simply provide a broader perspective in quadratic BSDEs, but also, by setting $f(|\cdot|)=0$, (partially) generalize  \cite{P1999}, \cite{BC2000}, \cite{B2003}, \cite{BLS2007}, etc.   Hence our approach can be seen as unified to the study of both quadratic BSDEs and non-quadratic BSDEs. 
Finally, as an application, we prove a probabilistic representation for the viscosity solution of the quadratic PDEs associated with the BSDEs of our study.

This paper is organized as follows.  In Section \ref{I}, we introduce some  functions used to treat the quadratic generator in (\ref{lq0}).
 In Section \ref{itokrylov}, we 
 prove  the It\^{o}-Krylov formula  and a generalized It\^{o} formula for $y\mapsto |y|^p (p\geq 1)$. The former one is used to treat discontinuous quadratic generators or  discontinuous quadratic growth,  and the later one is used to deduce the a priori estimates. 
  Section \ref{degenerate} reviews  purely quadratic BSDEs and studies their natural extensions, based on
 Bahlali et al \cite{BEK2013}.
 Section \ref{lpsolution}  concerns existence, comparison theorem, uniqueness, etc.
  Finally, in Section \ref{qpde}, we derive the nonlinear Feynman-Kac formula in our framework. 

Let us close this section by introducing all required notations.
  We fix  the time horizon $0<T<+\infty$ and   a $d$-dimensional Brownian motion  $(W_t)_{0\leq t\leq T}$ defined on a complete probability space $(\Omega, \mathcal{F},\mathbb{P})$.  $(\mathcal{F}_t)_{0\leq t \leq T}$ is the  filtration generated by $W$  and augmented by $\mathbb{P}$-null sets of $\mathcal{F}$. Any measurability will refer to this filtration. In particular, $\Prog$  denotes the progressive $\sigma$-algebra  on $\Omega \times [0, T]$. 
  
  As mentioned before, we only deal with $\mathbb{R}$-valued BSDEs of type (\ref{bsde}).  
We call the 
$\Prog\otimes\mathcal{B}(\mathbb{R})\otimes\mathcal{B}(\mathbb{R}^d)$-measurable random function  $F: \Omega\times[0, T]\times \mathbb{R}\times\mathbb{R}^d \rightarrow \mathbb{R}$ the \emph{generator} and the $\mathcal{F}_T$-measurable random variable $\xi$ the \emph{terminal value}.  
The conditions imposed on the generator are called the \emph{structure conditions}.
For notational convenience,  we sometimes write $(F, \xi)$ instead of (\ref{bsde}) to denote the BSDE with generator $F$ and terminal value $\xi$.

$\int_0^\cdot Z_sdW_s$, sometimes denoted by $Z\cdot W$,  refers to the vector stochastic integral; see, e.g.,
Shiryaev and  Cherny \cite{SC2002}. 
We call a process  $(Y, Z)$ valued in $\mathbb{R}\times\mathbb{R}^d$
a \emph{solution} of   (\ref{bsde}), if  $Y$ is a continuous adapted process and $Z$ is a $\Prog$-measurable process such that  $\mathbb{P}$-a.s.
$
\int_0^T |Z_s|^2 ds <+\infty$ and  $\int_0^T |F(s, Y_s, Z_s)| ds <+\infty,
$
and (\ref{bsde}) holds  $\mathbb{P}$-a.s. for any   $t\in[0, T]$.
The first inequality above ensures  that  $Z$ is integrable with respect to $W$ in the sense of vector stochastic integration.
 As a result, $Z\cdot W$ is a continuous local martingale.

 As will be seen later,  
 the BSDEs (\ref{bsde}) satisfying (\ref{lq0}) is solvable if $f(|\cdot|)$ belongs to $\mathcal{I}$, the set of integrable functions from $\mathbb{R}$ to $\mathbb{R}$
 which are bounded on any compact subset of $\mathbb{R}$.

For any random variable or process $Y$, we say $Y$ has some property if this is true except on a $\mathbb{P}$-null subset of $\Omega$. Hence we 
omit ``$\mathbb{P}$-a.s.'' in situations without ambiguity.  Define $\sgn(x): = \mathbb{I}_{\{x\neq 0 \}}\frac{x}{|x|}$.
For any c\`{a}dl\`{a}g adapted process $Y$, set $Y_{s, t}: = Y_t -Y_s$ and $Y^* : = \sup_{t\in [0, T]} |Y_t|$.
For any $\mathbb{R}$-valued $\Prog$-measurable process $H$, set $|H|_{s,t}:= \int_s^t H_u du$ and $|H|_t : = |H|_{0, t}$.   $\mathcal{T}$ stands for the set of  stopping times valued in $[0, T]$ and $\mathcal{S}$ denotes the space of continuous adapted processes. 
For any  local martingale $M$, we call $\{\sigma_n\}_{{n\in\mathbb{N}^+ }}\subset \mathcal{T}$ a \emph{localizing sequence} if  $\sigma_n$ increases stationarily to $T$ as $n$ goes to $+\infty$ and
$M_{\cdot \wedge \sigma_n}$ is a martingale for any $n\in\mathbb{N}^+$.
For later use,  we specify the following spaces under  $\mathbb{P}$.
\begin{itemize}
\item $\mathcal{S}^\infty$: the set of bounded processes in $\mathcal{S}$;
\item $\mathcal{S}^p (p\geq 1)$: the set of $Y\in \mathcal{S}$ with $Y^*\in \mathbb{L}^p$; 
\item $\mathcal{D}$: the set of $Y\in\mathcal{S}$ such that  $\{Y_\tau |  \tau\in\mathcal{T} \}$ is uniformly integrable;
\item $\mathcal{M}$: the space of $\mathbb{R}^d$-valued $\Prog$-measurable processes $Z$ such that  
$\mathbb{P}$-a.s.
$\int_0^T 
|Z_s|^2 ds <+\infty$;
for any $Z\in \mathcal{M}$, $Z\cdot W$ is a continuous local martingale;
\item $\mathcal{M}^p(p>0)$: the set of $Z\in\mathcal{M}$ with
\[
\norm{Z}_{\mathcal{M}^p}:=\mathbb{E}\Big[\Big(\int_0^T |Z_s|^2 ds \Big)^{\frac{p}{2}}\Big]^{\frac{1}{p}\wedge 1} <+\infty;
\]
\item $\mathcal{C}^p(\mathbb{R})$: the space of $p$ times continuously differentiable functions from $\mathbb{R}$ to $\mathbb{R}$;
\item$\mathcal{W}_{1, loc}^2(\mathbb{R})$: the Sobolev space  of measurable maps $u: \mathbb{R}\rightarrow \mathbb{R}$ such that both $u$ and its generalized derivatives $u^\prime, u^{\prime\prime}$ belong to $\mathbb{L}^1_{loc}(\mathbb{R})$.\\
\end{itemize}
 
 The above spaces 
are  Banach (respectively complete) under suitable norms (respectively metrics); we will not present these facts in more detail since
 they are not involved in our study. We call $(Y, Z)$ a $\mathbb{L}^p$ \emph{solution} of (\ref{bsde}) if $(Y, Z)$
 belongs to 
 $\mathcal{S}^p\times\mathcal{M}^p$.  This definition simply comes from the fact that the existence holds if $|\xi|+\int_0^T \alpha_s ds$ belongs to 
 $\mathbb{L}^p$. Analogously to most papers on $\mathbb{R}$-valued quadratic BSDEs, 
our existence result essentially relies  on 
 the monotone stability result of quadratic BSDEs; see, e.g.,  Kobylanski \cite{K2000} or  Briand and  Hu \cite{BH2008}.

\subsection{Functions of Class  $\mathcal{I}$}
\label{I}
In this section, we introduce the basic ingredients used to treat the quadratic generator in  (\ref{lq0}).
We recall that $\mathcal{I}$ is  the set of 
integrable  functions  from $\mathbb{R}$ to $\mathbb{R}$  which are bounded on any compact subset of $\mathbb{R}$.

{\bf $u^f$ Transform.} For any $f\in \mathcal{I}$, define   $u^f: \mathbb{R}\rightarrow \mathbb{R}$ and $M^f$ by
\begin{align*}
u^f (x)&: =  \int_0^x \exp{\Big(2\int_0^y f(u)du\Big)}dy,  \\
M^f&: = \exp\Big( 2\int_{-\infty}^{\infty} |f(u)| du\Big).
\end{align*}
Obviously, $1\leq M^f < + \infty$.
Moreover, the following properties hold by simple computations. Here we set $u:=u^f$ for notational convenience.
\begin{enumerate}
\item [\rm{(i)}] $u\in \mathcal{C}^1(\mathbb{R}) \cap \mathcal{W}_{1, loc}^2 (\mathbb{R})$ and 
$
u^{\prime\prime}(x) =2f(x)u^\prime(x)$  a.e.;
if  $f$ is continuous, then $u \in \mathcal{C}^2 (\mathbb{R})$;
\item [\rm{(ii)}] $u$ is strictly increasing and bijective from $\mathbb{R}$ to $\mathbb{R}$;
\item [\rm{(iii)}] $u^{-1} \in  \mathcal{C}^1(\mathbb{R}) \cap \mathcal{W}_{1, loc}^2 (\mathbb{R})$; if  $f$ is continuous, then $u^{-1} \in \mathcal{C}^{2}(\mathbb{R})$;
\item [\rm{(iv)}] 
$
\frac{|x|}{M}  \leq |u(x)| \leq M|x|
$ 
and
$\frac{1}{M}\leq u^\prime(x)\leq M.$
\end{enumerate}

{\bf $v^f$ Transform.} For any $f\in \mathcal{I}$, define $v^f: \mathbb{R}\rightarrow \mathbb{R}^+$ by
\begin{align*}
v^f(x): = \int_0^{|x|} u^{(-f)}(y) 
\exp\Big(2\int_0^y f(u)du     \Big)
  dy.
\end{align*}
Set $v:=v^f.$
Simple computations  give
\begin{enumerate}
\item [(i)]$v\in \mathcal{C}^1(\mathbb{R}) \cap \mathcal{W}_{1, loc}^2 (\mathbb{R})$  and
$v^{\prime\prime}(x)-2f(|x|)|v^\prime(x)|=1$ a.e.;
if  $f$ is continuous, then 
$v\in \mathcal{C}^2(\mathbb{R});$
\item [(ii)] $v(x)\geq 0, \sgn(v^\prime(x))=\sgn(x)$ and $v^{\prime\prime}(0)=1;$
 \item [(iii)]
 $
 \frac{x^2}{2M^2}\leq v(x) \leq \frac{M^2x^2}{2}$ and $\frac{|x|}{M^2}\leq|v^\prime(x)|\leq M^2 |x|$.
\end{enumerate}
In the sequel of our study, $u^f$ and $v^f$ exclusively stand for the above transforms associated with $f\in\mathcal{I}$. Hence in situations without ambiguity, we  denote $u^f, v^f, M^f$ by $u, v, M$, respectively.
\subsection{Krylov Estimate and  the It\^{o}-Krylov Formula}
\label{itokrylov}
The first auxiliary result  is
 the Krylov estimate.
Later, it is used to 
 prove an  It\^{o}'s-type formula for functions in $\mathcal{C}^1(\mathbb{R})\cap \mathcal{W}_{1, loc}^2(\mathbb{R})$. This helps to deal with 
 (possibly discontinuous) quadratic generators.
As the second application, we derive a generalized It\^{o} formula for $y\mapsto |y|^p(p\geq 1)$ which is not smooth enough for $1\leq p <2$. This is a basic tool to study $\mathbb{L}^p(p\geq 1)$ solutions.

 To allow the existence of a local time  in  particular situations, we  study equations of type
\begin{align}
Y_t = \xi + \int_t^T F(s, Y_s, Z_s) ds +\int_t^T dC_s -\int_t^T Z_s dW_s, \label{fc}
\end{align}
where $C$ is a continuous adapted process of finite variation. We denote its total variation process by $V_\cdot(C)$.  Likewise,  sometimes  we denote   (\ref{fc}) by 
$(F, C, \xi)$. The solution of (\ref{fc}) is defined analogously to that 	of   (\ref{bsde}).

Now we prove the Krylov estimate for (\ref{fc}).
A more complicated version  not needed  for our study can be found in Bahlali et al \cite{BEO2014}. 

\begin{Lemma}[Krylov Estimate] \label{krylovestimate} 
Consider {\rm(\ref{fc})}.
For any measurable function $\psi: \mathbb{R}\rightarrow \mathbb{R}^+$,
\begin{align}
\mathbb{E} \Big[  \int_0^{ \tau_{m}} \psi(Y_s)|Z_s|^2 ds \Big] \leq 6m \norm{\psi}_{\mathbb{L}^1([-m, m])}, \label{psifun}
\end{align}
where  $\tau_m$ is a stopping time defined by 
\begin{align*}
\tau_{m} : = \inf\Big \{t \geq 0: |Y_t| + V_t(C) + \int_0^t |F(s, Y_s, Z_s)| ds \geq m \Big\} \wedge T.
\end{align*}
\end{Lemma}
\begin{Proof}
Without loss of generality we assume
$\norm{\psi}_{\mathbb{L}^1([-m, m])} <+\infty$. 
For   each $n\in\mathbb{N}^+$, set 
\begin{align*}
\tau_{m,n} &: = \tau_{m} \wedge \inf \Big\{t\geq 0: \int_0^t |Z_s|^2 ds \geq n \Big\}.
\end{align*} 
Let $a\in [-m, m]$. By Tanaka's formula, 
\begin{align}
(Y_{t\wedge \tau_{m,n}} - a)^{-}  & =  (Y_0 - a)^- - \int_0^{t\wedge\tau_{m,n}} \mathbb{I}_{\{Y_s < a\} }dY_s + \frac{1}{2}L_{t\wedge \tau_{m,n}}^a(Y)\nonumber\\ 
& =  (Y_0 - a)^- + \int_0^{t\wedge \tau_{m,n}} \mathbb{I}_{\{Y_s < a\}} F(s, Y_s, Z_s)ds +\int_0^{t\wedge \tau_{m, n}} \mathbb{I}_{\{ Y_s < a\}}dC_s \nonumber\\
&-\int_0^{t\wedge\tau_{m,n}}\mathbb{I}_{\{Y_s < a\}} Z_s dW_s 
+\frac{1}{2} L_{{t\wedge \tau_{m,n}}}^a (Y), \label{ito-}
\end{align}
where $L^{a}(Y)$ is the local time of $Y$ at $a$.
To estimate the local time, we put it on the left-hand side and the rest terms on the right-hand side. 
Since $x \mapsto (x-a)^-$ is Lipschitz-continuous,  we deduce from the definition of $\tau_{m, n}$ that
\[
(Y_0 - a)^- - (Y_{t\wedge \tau_{m, n}} - a)^- \leq | Y_0 -Y_{t\wedge{\tau_{m, n}}}| \leq 2m.
\]
Meanwhile, the definition of $\tau_m$ also implies that the sum of the $ds$-integral and $dC$-integral is bounded by $m$. Hence, we have
\[
\mathbb{E}\big[L_{t\wedge \tau_{m,n}}^a(Y)\big] \leq 6m. 
\]
By Fatou's lemma applied to the sequence indexed by $n$, 
\[
\sup_{a\in[-m,  m]}\mathbb{E}\big[L_{t\wedge \tau_{m}}^a(Y)\big] \leq 6m.
\]
We then use  
time occupation formula for continuous semimartingales (see Chapter VI., Revuz and  Yor  \cite{RY1999})  and the above inequality to obtain
\begin{align*}
\mathbb{E}\Big[ \int_0^{T\wedge\tau_{m}} \psi(Y_s)|Z_s|^2 ds  \Big] 
 &= 
\mathbb{E}\Big[ \int_{-m}^{m} \psi(x)L_{T\wedge \tau_m}^x(Y)dx  \Big]\\
&= \int_{-m}^m \psi(x) \mathbb{E}\big[L_{T\wedge \tau_{m}}^x (Y)\big] dx\\
&\leq 6m \norm{\psi}_{\mathbb{L}^1([-m, m])}.
\end{align*}
\qed
\end{Proof}

As an immediate consequence of Lemma 
\ref{krylovestimate}, we have $\mathbb{P}$-a.s.
\begin{align}
\int_0^T \mathbb{I}_{\{ Y_s \in A \}}|Z_s|^2 ds = 0, \label{Aneg}
\end{align}
 for any $A\subset\mathbb{R}$ with null Lebesgue measure.  This will be used later several times.

  Given Lemma \ref{krylovestimate}, we turn to the main results of this section. The following generalized It\^o formula is proved in Bahlali et al \cite{BEO2014}.

\begin{thm}[It\^{o}-Krylov Formula]\label{ik} If $(Y, Z)$ is a solution of {\rm(\ref{fc})}, then
 for any  $u\in\mathcal{C}^1(\mathbb{R})\cap \mathcal{W}_{1, loc}^2(\mathbb{R})$,  we have $\mathbb{P}$-a.s. for all $t\in [0, T]$, 
\begin{align}
u(Y_t) = u(Y_0) + \int_0^t u^\prime (Y_s) dY_s +\frac{1}{2} \int_0^t u^{\prime\prime}(Y_s)|Z_s|^2 ds. \label{ikformula} 
\end{align}

\end{thm}
\begin{Proof}
We use $\tau_{m}$  defined in Lemma \ref{krylovestimate} (Krylov estimate). 
Note that  $\tau_{m}$ increases stationarily to $T$ as $m$ goes to $+\infty$. It is therefore sufficient to prove the equality for $u(Y_{t\wedge\tau_m})$. To this end we use an approximation procedure. We consider $m$ such that  $\mathbb{P}$-a.s. $m\geq |Y_0|$.
Let $u_n$ be a sequence of functions in  $\mathcal{C}^2(\mathbb{R})$ satisfying
\begin{enumerate}
\item [(i)] $u_n$ converges uniformly to $u$ on $[-m, m]$; 
\item [(ii)] $u_n^\prime$ converges uniformly to $u^\prime$ on $[-m, m]$;
\item [(iii)]$u_n^{\prime\prime}$ converges in $\mathbb{L}^1([-m, m])$ to $u^{\prime\prime}$.
\end{enumerate}
By It\^{o}'s formula,
\[
u_n(Y_{t\wedge \tau_{m}}) = u_n(Y_0) + \int_0^{t\wedge \tau_{m}} u_n^\prime (Y_s) dY_s +\frac{1}{2} \int_0^{t\wedge \tau_{m}} u_n^{\prime\prime}(Y_s)|Z_s|^2 ds.
\]
Due to  (i) and $|Y_{t \wedge \tau_m}| \leq m$, 
$u_n (Y_{\cdot \wedge \tau_{m}})$  converges to $u(Y_{\cdot \wedge \tau_{m}})$ $\mathbb{P}$-a.s.  
uniformly on $[0, T]$
as $n$ goes to $+\infty$; the second term converges in probability to 
\[
\int_0^{t \wedge \tau_{m}} u^\prime (Y_s) dY_s
\]
by (ii)
 and dominated convergence for stochastic integrals; the last term converges in probability to 
\[
\frac{1}{2}\int_0^{t\wedge \tau_{m}} u^{\prime\prime}(Y_s) |Z_s|^2 ds
\]
due to  (iii) and Lemma \ref{krylovestimate}. 
Indeed,  Lemma \ref{krylovestimate} implies 
\begin{align*}
\mathbb{E}\Big[\int_0^{\tau_m} |u^{\prime\prime}_n-u^{\prime\prime}|(Y_s)|Z_s|^2 ds  \Big]
\leq 6m \norm{u^{\prime\prime}_n-u^{\prime\prime}}_{{\mathbb{L}^1{([-m, m])}}}.
\end{align*} Hence  collecting these convergence results  gives (\ref{ikformula}). 
By the continuity of both sides of (\ref{ikformula}), the quality also holds $\mathbb{P}$-a.s. for all $t\in[0, T]$.
\qed
\end{Proof}

To study $\mathbb{L}^p(p\geq 1)$ solutions we now prove an It\^o's-type formula for $y \mapsto |y|^p(p\geq 1)$ which
is not smooth enough for $1\leq p<2$.  The proof for 
 multidimensional It\^{o} processes can be found, e.g., in  Briand et al \cite{B2003}. In contrast to their approach, 
we give a novel and simpler proof  for BSDE framework but point out that it can be also extended to  It\^{o} processes. 
\begin{Lemma}\label{lpinequality}
Let $p \geq 1$. If $(Y, Z)$ is a solution of {\rm(\ref{fc})},
then we have $\mathbb{P}$-a.s. for all $t\in [0, T]$, 
\begin{align}
|Y_t|^p &+ \frac{p(p-1)}{2} \int_t^T  \mathbb{I}_{\{ Y_s \neq 0\}}|Y_s|^{p-2}  |Z_s|^2ds  \nonumber\\
&   = |\xi|^p - p \int_t^T \sgn(Y_s)|Y_s|^{p-1} dY_s-\mathbb{I}_{\{ p =1 \}}\int_t^T dL_s^0(Y), \label{lemma23}
\end{align}
where $L^0(Y)$ is the local time of 
$Y$ at $0$. 
\end{Lemma}
\begin{Proof}
(i). $p = 1 $. This is immediate from Tanaka's formula. 

(ii). $ p> 2$. $y \mapsto |y|^p \in \mathcal{C}^2(\mathbb{R})$. Hence this is immediate from It\^{o}'s formula.

(iii). $p=2$.  $y \mapsto |y|^p \in \mathcal{C}^2(\mathbb{R})$. 
Due to (\ref{Aneg}), 
$\int_0^\cdot|Y_s|^{p-2}  |Z_s|^2ds$
is indistinguishable from $\int_0^\cdot  \mathbb{I}_{\{ Y_s \neq 0\}}|Y_s|^{p-2}  |Z_s|^2ds$. By taking this fact into account, 
this equality is thus immediate from It\^o's formula.
    
(iv). $1<p<2$. We use an approximation argument. Define 
\[
u_\epsilon (y): = \big(y^2 + \epsilon^2\big)^{\frac{1}{2}}.
\]
Hence for any $\epsilon >0$, we have $u_\epsilon^p\in \mathcal{C}^2(\mathbb{R}).$ By It\^{o}'s formula, 
 \begin{align}
 u^p_\epsilon (Y_t) &=  u^p_\epsilon (\xi) - p \int_t^T Y_s u^{p-2}_\epsilon (Y_s) d Y_s -   \frac{1}{2} \int_t^T \big( p u_\epsilon^{p-2}(Y_s)+p(p-2)|Y_s|^2u_\epsilon^{p-4}(Y_s)\big)
   |Z_s|^2   ds.  \label{uepsilon}
 \end{align} Now we send $\epsilon$ to $0$.
 $u_\epsilon(y) \longrightarrow |y|$ pointwise implies
 $
u_\epsilon(Y_t)^p \longrightarrow
|Y_t|^p$ and  $\ u_\epsilon(\xi)^p
\longrightarrow |\xi|^p$ pointwise on $\Omega$.
  Secondly, 
$yu_\epsilon^{p-2}(y) \longrightarrow \sgn(y)|y|^{p-1}$ pointwise implies 
by dominated convergence for stochastic integrals that
\begin{align*}
&\int_t^T Y_s \sgn(Y_s) u^{p-2}_\epsilon (Y_s) dY_s {\longrightarrow}\int_t^T |Y_s|^{p-1}dY_s \ \text{in probability.}
\end{align*}
To prove that the $ds$-integral in (\ref{uepsilon}) also converges, we split it into two parts and argue their convergence respectively.
 Note that  
 \begin{align}
pu_\epsilon^{p-2}(Y_s)+p(p-2)|Y_s|^2u_\epsilon^{p-4}(Y_s) = p\epsilon^2 u_\epsilon^{p-4}(Y_s) + p(p-1) |Y_s|^2 u_\epsilon^{p-4} (Y_s). \label{2ndorder}
\end{align}
For the second term on the right-hand side of (\ref{2ndorder}), we have
\[
|Y_s|^2 u_\epsilon^{p-4} (Y_s)  =  \mathbb{I}_{\{Y_s \neq 0\}}  |Y_s|^{p-2} \Big|\frac{|Y_s|}{u_\epsilon(Y_s)} \Big|^{4-p}.
\]
Since
$
\frac{|y|}{u_\epsilon(y)} {\nearrow} \mathbb{I}_{\{ y \neq 0\}}
$, 
 monotone convergence gives
\[
\int_t^T |Y_s|^2 u_\epsilon^{p-4}(Y_s) |Z_s|^2ds
\longrightarrow 
\int_t^T \mathbb{I}_{\{Y_s \neq 0 \}}|Y_s|^{p-2}|Z_s|^2  ds \  \text{pointwise on}\ \Omega. 
\]
It thus remains to prove  the $ds$-integral concerning the first term on the right-hand side of (\ref{2ndorder})
  converges to $0$. To this end,  
 we use Lemma \ref{krylovestimate} (Krylov estimate)
 and the same localization procedure.
This gives
\begin{align*}
\mathbb{E} \Big[ \int_0^{\tau_m} \epsilon^2 u_\epsilon^{p-4}(Y_s) |Z_s|^2ds  \Big] &\leq 6m \epsilon^2 \int_{-m}^{m} (x^2 + \epsilon^2)^{\frac{p-4}{2}} dx\\
& \leq  12m\epsilon^2\int_0^m  (x^2 + \epsilon^2)^{\frac{p-4}{2}} dx\\
& \leq  12\cdot 2^{\frac{4-p}{2}}m\epsilon^2 \int_{0}^{m} (x+\epsilon)^{p-4} dx\\
& \leq  12\cdot 2^{\frac{4-p}{2}}m\epsilon^2 \int_{\epsilon}^{m+\epsilon} x^{p-4} dx\\
&= \frac{12\cdot 2^{\frac{4-p}{2}} m}{p-3} \big (\epsilon^2  (m+\epsilon)^{p-3} - \epsilon^{p-1}\big),
\end{align*}
which, due to $1<p$, converges to $0$ as $\epsilon $ goes to $0$. Hence 
$
\int_0^{\cdot}\epsilon^2 u_\epsilon^{p-4}(Y_s)|Z_s|^2 ds 
$
 converges $u.c.p$ to $0$. 
Collecting all convergence results above gives (\ref{lemma23}). Finally, the continuity
 of each term in (\ref{lemma23}) implies that the equality
also holds $\mathbb{P}$-a.s. for all $t\in [0, T]$.
\qed
\end{Proof}
\subsection{$\mathbb{L}^p(p \geq 1)$ Solutions of Purely Quadratic BSDEs}
\label{degenerate}
Before turning to the main results of this paper, we partially extend the existence and uniqueness result for
purely quadratic BSDEs studied by
Bahlali et al \cite{BEO2014}. Later,  we present their natural extensions and the motivations of our work. 
These BSDEs  are called purely quadratic, since the generator takes the form
$F(t, y, z)=f(y)|z|^2$.  
The solvability simply comes from the function $u^f$ defined in Section \ref{I} which transforms  better known BSDEs to  $(f(y)|z|^2, \xi)$ by It\^{o}-Krylov formula. 
\begin{thm}
\label{degeneratethm}
Let $f\in \mathcal{I}$ and $\xi \in \mathbb{L}^p (p \geq 1)$. Then there exists a unique solution of
\begin{align}
Y_t = \xi + \int_t^T f(Y_s)|Z_s|^2ds -\int_t^T Z_s dW_s \label{pqbsde}.
\end{align}
Moreover, if $p>1$, the solution belongs to $\mathcal{S}^p\times\mathcal{M}^p$; if $p=1$, the solution belongs to $\mathcal{D}\times\mathcal{M}^q$ for any $q\in(0, 1)$.
\end{thm}
\begin{Proof}
Let $u: = u^f$ and $M:=M^f.$
Then $u, u^{-1}\in
\mathcal{C}^1(\mathbb{R})\cap \mathcal{W}_{1, loc}^2(\mathbb{R})$. The existence and uniqueness result can be seen as a one-on-one correspondence between solutions of  BSDEs.

(i). Existence. $|u(x)| \leq M |x|$ implies  
$u(\xi) \in \mathbb{L}^p$.
By  It\^{o} representation theorem,  
there exists a unique pair $(\widetilde{Y}, \widetilde{Z})$ which solves $(0, u(\xi))$, i.e.,
\begin{align}
d\widetilde{Y}_t = \widetilde{Z}_t dW_t, \ \widetilde{Y}_T = u(\xi). \label{tbsde}
\end{align}
We aim at proving
\begin{align}(Y, Z):=
(u^{-1}(\widetilde{Y}),
\frac{\widetilde{Z}}{u^\prime(u^{-1}(\widetilde{Y}))}) \label{yz}
\end{align}
solves (\ref{pqbsde}).  
 It\^o-Krylov formula applied to $Y_t = u^{-1}(\widetilde{Y}_t)$ yields 
\begin{align}
dY_t =\frac{1}{u^\prime(u^{-1}(\widetilde{Y}_t))} d\widetilde{Y}_t - \frac{1}{2}\Big(\frac{1}{u^\prime(u^{-1}(\widetilde{Y}_t))}\Big)^2 \frac{u^{\prime\prime}(u^{-1}(\widetilde{Y}_t))}{u^\prime(u^{-1}(\widetilde{Y}_t))} |\widetilde{Z}_s|^2 ds. \label{ytt}
\end{align}
To simplify (\ref{ytt}) let us recall that $u^{\prime\prime}(x)= 2f(x)u^{\prime}(x)$ a.e. Hence (\ref{yz}), (\ref{ytt}) and (\ref{Aneg}) give
\[
dY_t = -f(Y_t)|Z_t|^2 dt +  Z_t dW_t, \ Y_T = \xi,
\]
i.e., $(Y, Z)$ solves (\ref{pqbsde}).

(ii). Uniqueness. Suppose $(Y, Z)$ and  $(Y^\prime, Z^\prime)$ are solutions of (\ref{pqbsde}). By It\^o-Krylov formula applied to $u(Y)$ and  $u(Y^\prime)$, we deduce that 
$(u(Y), u^\prime(Y)Z)$ and $(u(Y^\prime), u^\prime(Y^\prime)Z^\prime)$
 solve
$(0, u(\xi))$. But from  (i) it is known that they coincide. Transforming $u({Y})$ and $u(Y^\prime)$ via the bijective function $u^{-1}$ yields the uniqueness result.

(iii). We prove the estimate for the unique solution $(Y, Z)$. For $p>1$, Doob's $\mathbb{L}^p(p>1)$ maximal inequality used to (\ref{tbsde}) implies  $(\widetilde{Y}, \widetilde{Z})\in\mathcal{S}^p\times\mathcal{M}^p$. Hence $({Y}, {Z})\in\mathcal{S}^p\times\mathcal{M}^p$, due to $|u^\prime(x)|\geq \frac{1}{M}$ and $|u^{-1}(x)| \leq M |x|$.
 For $p=1$, $\widetilde{Y}\in \mathcal{D}$ since it is a martingale on $[0, T]$. In view of the above properties of $u$ we have $Y\in\mathcal{D}$. The estimate for $Z$ is immediate from
Lemma 6.1, Briand et al \cite{B2003} which is a version of $\mathbb{L}^p(0<p<1)$ maximal inequality for martingales.
\qed
\end{Proof}
\begin{remm}
If $\xi$ is a general $\mathcal{F}_T$-measurable random variable,
 Dudley representation theorem (see Dudley \cite{d1977}) implies that 
 there  still exists a solution of (\ref{tbsde}) and hence a solution of (\ref{pqbsde}). However, the solution in general is not unique.
 
The proof of Theorem \ref{degeneratethm} indicates that 
$f$  being bounded on compact subsets of $\mathbb{R}$ is not needed 
for the existence and uniqueness result of purely quadratic BSDEs.
\end{remm}
\begin{prop}[Comparison]
\label{dcompare}
Let $f, g\in \mathcal{I}$, $\xi, \xi^\prime \in\mathbb{L}^p(p\geq 1)$ and
 $(Y, Z)$, $(Y^\prime, Z^\prime)$ be the unique solutions of $(f(y)|z|^2, \xi)$, $(g(y)|z|^2, \xi^\prime)$, respectively. 
If $f\leq g$ a.e. and $\mathbb{P}$-a.s. $\xi\leq \xi^\prime$, then  $\mathbb{P}$-a.s. $Y_\cdot \leq Y^\prime_\cdot$. 
\begin{Proof}
Again we transform so as to compare better known BSDEs. Set $u:=u^f$. For any $\tau\in\mathcal{T}$, 
 It\^{o}-Krylov formula yields
\begin{align*}
u (Y_{t\wedge \tau}^\prime)  
&=
u(Y_{\tau}^\prime) + \int_{t\wedge \tau}^\tau \Big(u^\prime(Y_s^\prime)g(Y_s^\prime) |Z_s^\prime|^2
-\frac{1}{2}u^{\prime\prime}(Y_s^\prime)|Z_s^\prime|^2\Big)
 ds 
-\int_{t\wedge\tau}^\tau u^\prime (Y_s)Z_s^\prime dW_s.
\\
&=
u(Y_{\tau}^\prime)+\int_{t\wedge \tau}^\tau u^\prime (Y_s^\prime) 
\big(g(Y_s^\prime)-f(Y_s^\prime)\big) |Z_s^\prime|^2 ds
-\int_{t\wedge \tau}^\tau u^\prime (Y_s)Z_s^\prime dW_s \\
&\geq
u(Y_{\tau}^\prime) -\int_{t\wedge \tau}^\tau u^\prime (Y_s)Z_s^\prime dW_s,
\end{align*}
where the last two lines are due to
$u^{\prime\prime}(x)  = 2 f(x)u^\prime(x)$ a.e., $g\geq f$ a.e. and (\ref{Aneg}). In the next step, we want to eliminate the local martingale part by a localization procedure. Note that 
$\int_t^\cdot u^{\prime}(Y_s)Z_s^\prime dW_s$ is a local martingale on $[t, T]$. Set 
$\{\tau_n\}_{n\in\mathbb{N}^+}$
to be its localizing sequence on $[t, T]$. Replacing $\tau$ by $\tau_n$ in the above inequality thus gives $\mathbb{P}$-a.s.
\[
u(Y_t^\prime) \geq \mathbb{E}\big[u(Y_{t\wedge \tau_n}^\prime)\big|\mathcal{F}_t\big].
\]
This implies that, for any $A\in\mathcal{F}_t$,  we have
\[
\mathbb{E}\big[ u(Y_t^\prime)\mathbb{I}_A\big] \geq 
\mathbb{E}\big[ u(Y_{t\wedge\tau_n}^\prime)\mathbb{I}_A  \big].
\]
Since $u(Y^\prime)\in\mathcal{D}$, we can use Vitali convergence theorem to obtain
\[
\mathbb{E}\big[ u(Y_t^\prime)\mathbb{I}_A \big]
\geq \mathbb{E}\big[ u(\xi^\prime)  \mathbb{I}_A\big] =
\mathbb{E}\big[\mathbb{E}\big[ u(\xi^\prime)\big| \mathcal{F}_t \big]  \mathbb{I}_A\big].
\]
Note that this inequality holds for any $A\in\mathcal{F}_t$. Hence, by choosing $A=\{u(Y_t^\prime) < \mathbb{E}[u(\xi^\prime)|\mathcal{F}_t]\}$, we obtain
$
u(Y_t^\prime) \geq \mathbb{E}\big[u(\xi^\prime)\big|\mathcal{F}_t\big].$  Since $\xi^\prime \geq \xi$ and $u$ is increasing, we further have 
$u(Y_t^\prime) \geq 
\mathbb{E}\big[ u(\xi)\big|\mathcal{F}_t\big].$ Let us recall that, by Theorem \ref{degeneratethm}, 
$(u(Y), u^\prime (Y)Z)$ is the unique solution of $(0, u(\xi))$. Hence, $
u(Y_t^\prime) \geq u(Y_t).$
Transforming both sides via the bijective increasing function $u^{-1}$
yields  $\mathbb{P}$-a.s. $Y_t\leq Y_t^\prime$. By the continuity of $Y$ and $Y^\prime$ we have $\mathbb{P}$-a.s. $Y_\cdot\leq Y^\prime_\cdot$.
\qed
\end{Proof}
\end{prop}
\begin{remm}
In   Proposition \ref{dcompare}, we rely on the fact that 
$\mathbb{P}$-a.s.
\begin{align}\int_0^\cdot 
\Big(\frac{1}{2}u^{\prime\prime}(Y_s^\prime) - f(Y_s^\prime) u^{\prime}(Y_s^\prime)\Big)|Z_s^\prime|^2 ds=0, \label{aac}
\end{align}
 even though   
$u^{\prime\prime}(x)  = 2 f(x)u^\prime(x)$  only holds  almost everywhere on $\mathbb{R}$.
Here we prove it. 
 Let  $A$ be the subset of $\mathbb{R}$ on which $u^{\prime\prime}(x)  = 2 f(x)u^\prime(x)$ fails.  Hence, 
\[\int_0^\cdot \mathbb{I}_{\{ Y^\prime_s\in \mathbb{R}\backslash A\}}
\Big|\frac{1}{2}u^{\prime\prime}(Y_s^\prime) - f(Y_s^\prime) u^{\prime}(Y_s^\prime)\Big||Z_s^\prime|^2 ds=0.\]
Meanwhile, by (\ref{Aneg}), we have 
$\mathbb{P}$-a.s.
\[\int_0^\cdot \mathbb{I}_{\{ Y^\prime_s\in A\}}
\Big|\frac{1}{2}u^{\prime\prime}(Y_s^\prime) - f(Y_s^\prime) u^{\prime}(Y_s^\prime)\Big||Z_s^\prime|^2 ds=0.\]
Hence, (\ref{aac}) holds 
$\mathbb{P}$-a.s.
 This fact also applies to Theorem \ref{degeneratethm} and all results in the sequel of our study.  
\end{remm}

Theorem \ref{degeneratethm}   and Proposition \ref{dcompare} are
based on a one-on-one correspondence between solutions (respectively the unique solution) of BSDEs. Hence it is natural to generalize as follows.  Set  $f\in\mathcal{I}, u:=u^f$,  $F(t, y, z):=  G(t, y, z) + f(y)|z|^2$  and
\begin{align}
\widetilde{F} (t, y, z): =
u^\prime (u^{-1}(y)) G(t, u^{-1}(y),\frac{z}{u^\prime (u^{-1}(y))}).\label{adds}
\end{align}
If $G$ ensures the existence of a solution of   $(\widetilde{F}, u(\xi))$, we can transform it via ${u}^{-1}$
to  a solution of $(F, \xi)$. An example is that
$G$ is of continuous linear growth in $(y, z)$ where 
the existence of a maximal (respectively minimal) solution of $(\widetilde{F}, u(\xi))$ can be proved in the spirit of Lepeltier and San Martin \cite{LS1997}.

When the generator is continuous in $(y, z)$, a more general situation is    linear-quadratic growth, i.e., 
\begin{align}
|H(t, y, z)| \leq \alpha + \beta|y| + \gamma |z| + f(|y|)|z|^2 : =F(t, y, z), \label{lq}
\end{align}
for some $\alpha, \beta, \gamma \geq 0$. 
The existence result  then consists of viewing the maximal (respectively minimal) solution of $(F, \xi^+)$ (respectively $(-F, -\xi^-)$) as the a priori bounds for solutions of $(H, \xi)$, 
and using a combination of a localization procedure and  the monotone stability result developed by  Briand and  Hu \cite{BH2006}, \cite{BH2008}.
For details the reader shall refer to Bahlali et al \cite{BEO2014}.

However, either an additive structure in (\ref{adds}) or a linear-quadratic growth (\ref{lq})
 is too restrictive 
   and  uniqueness  is not available in general.  Considering this limitation, we devote Section \ref{lpsolution} to  the solvability under  milder structure conditions.  
\subsection{$\mathbb{L}^p(p > 1)$ Solutions of Quadratic BSDEs}
\label{lpsolution} 
With the preparatory work in Section \ref{I}, \ref{itokrylov}, \ref{degenerate}, we study  $\mathbb{L}^p (p > 1)$ solutions of quadratic BSDEs under general assumptions. We deal with the quadratic generators in the spirit of   
 Bahlali et al	\cite{BEO2014}, 
derive the a priori estimates in the spirit of Briand et al \cite{B2003} and prove the existence and uniqueness result in the spirit of  Briand et al \cite{BH2006}, \cite{BH2008}, \cite{BLS2007}. This section can also be seen as a generalization of these works.  
The following assumptions on $(F, \xi)$  ensure the estimates and an existence result. 
\begin{as}\label{lpas1}
Let $p\geq 1$.
There exist  $\beta \in \mathbb{R}, \gamma \geq 0$, an $\mathbb{R}^+$-valued $\Prog$-measurable process $\alpha$, $f(|\cdot|)\in\mathcal{I}$  and a continuous nondecreasing function $\varphi: \mathbb{R}^+ \rightarrow \mathbb{R}^+$ with $\varphi(0)=0$ such that     
$  |\xi|  + |\alpha |_T \in \mathbb{L}^p$ and   $\mathbb{P}$-a.s.
  \begin{enumerate} 
\item [(i)] for any $t\in[0, T]$, $(y, z)\longmapsto F(t, y, z)$ is continuous; 
  \item [(ii)] $F$ is ``monotonic''  at  $y=0$, i.e.,  for any $(t,y, z)\in [0, T]\times\mathbb{R}\times\mathbb{R}^d$,
    \[\sgn (y) F(t, y, z) \leq  \alpha_t + \beta |y| + \gamma |z| + f(|y|) |z|^2; \]
  \item [(iii)]  for any $(t, y, z)\in[0, T]\times\mathbb{R}\times\mathbb{R}^d$,
  \[|F(t, y, z)| \leq \alpha_t + \varphi(|y|) + \gamma |z| + f(|y|)|z|^2.\] 
\end{enumerate}
\end{as} 

It
is worth noticing that,
given \ref{lpas1}(iii) and $f(|\cdot|)=0$, 
\ref{lpas1}(ii) is a consequence of  $F$ being monotonic at $y=0$. Indeed,  
\[
\sgn (y-0)\big( F(t, y, z) - F(t, 0, z)\big) \leq \beta|y| 
\]
implies
\begin{align*}
\sgn(y) F(t, y, z) &\leq
F(t, 0, z) + \beta|y| \\
&\leq \alpha_t + \beta |y| +\gamma |z|.
\end{align*}
This explains why  we keep  saying that  $F$ is monotonic  at $y=0$, even though $y$ also appears in the quadratic term.
Secondly, our results don't rely on the specific choice of $\varphi$. 
Hence the growth condition in $y$ can be arbitrary as long as \ref{lpas1}(i)(ii) hold.
Assumptions of this type for different settings can also be found in, e.g.,   \cite{BC2000}, \cite{B2003}, \cite{BH2008}.
 Finally,  $f(|\cdot|)$ can be discontinuous; $f(|\cdot|)$ being $\mathbb{R}^+$-valued appears more naturally in the growth condition.
\begin{Lemma}
[A Priori Estimate (i)]
\label{lpestimate1} Let $p \geq 1$ and  {\rm \ref{lpas1}} hold for $(F, \xi)$. If $(Y, Z) \in \mathcal{S}^p \times \mathcal{M} $ is a solution of $(F, \xi)$, then 
\begin{align*}
\mathbb{E}\Big[\Big( \int_0^T |Z_s|^2 ds \Big)^{\frac{p}{2}}\Big] + \mathbb{E}\Big[ \Big(\int_0^T  f(|Y_s|) |Z_s|^2 ds \Big)^{p} \Big] 
\leq c\Big( \mathbb{E}\big[(Y^*)^p + |\alpha|_T^p\big]\Big),
\end{align*}
where  $c$ is a  constant  only depending on
$T, M^{f(|\cdot|)}, \beta, \gamma, p$. 
\end{Lemma}
\begin{Proof} 
Set $v:=v^{f(|\cdot|)}$ and $M:=M^{f(|\cdot|)}$.
For any $\tau\in\mathcal{T}$,  It\^o-Krylov formula  yields
\begin{align}
v(Y_0) &= v(Y_\tau)
+\int_0^\tau  v^\prime(Y_s) F(s, Y_s, Z_s) ds \nonumber
\\ 
&-\frac{1}{2}\int_0^\tau v^{\prime\prime}(Y_s)|Z_s|^2 ds 
- \int_0^\tau v^\prime(Y_s)Z_s dW_s.\label{itov1}
\end{align}
Due to $\sgn(v^\prime(x))=\sgn(x)$ and \ref{lpas1}(ii),  we have
\begin{align}
v^\prime(Y_s)F(s, Y_s, Z_s)
\leq 
|v^\prime(Y_s)|\big(\alpha_t +\beta|Y_s|+\gamma |Z_s|+ f(|Y_s|)|Z_s|^2\big). \label{itov0}
\end{align}
Recall that
$v^{\prime\prime}(x)-2f(|x|)|v^\prime(x)|=1$ a.e.
Hence (\ref{itov1}) and (\ref{itov0}) 
give
\begin{align*}
\frac{1}{2}\int_0^\tau |Z_s|^2 ds \leq v(Y_\tau)
+\int_0^\tau |v^\prime(Y_s)|
\big(\alpha_s +\beta|Y_s| +\gamma |Z_s| \big) ds -\int_0^\tau v^\prime(Y_s)Z_s dW_s.
\end{align*}
Moreover, since $v(x)\leq \frac{M^2x^2}{2}$ and  $|v^\prime(x)|\leq M^2|x|$,   this inequality gives
\begin{align}
\int_0^\tau |Z_s|^2 ds \leq c_1 (Y^*)^2
+c_1\int_0^\tau  |Y_s|
\big(\alpha_s +|Y_s| + |Z_s| \big) ds -2\int_0^\tau v^\prime(Y_s)Z_s dW_s, \label{itov2}
\end{align}
where $c_1:= 2{M^2}(1 \vee \beta\vee \gamma)$.
Note that  in (\ref{itov2}),
\begin{align*}
\int_0^{\tau} |Y_s| \alpha_s ds &\leq   \frac{1}{2}(Y^*)^2 + \frac{1}{2}|\alpha|_{T}^2, \\
c_1 \int_0^{\tau} |Y_s| |Z_s| ds &\leq  \frac{1}{2}c^2_1 T\cdot (Y^*)^2 + \frac{1}{2}\int_0^{\tau} |Z_s|^2 ds.
\end{align*}
Hence (\ref{itov2}) yields
\[
\int_0^{\tau}  |Z_s|^2 ds  \leq (3c_1+c^2_1 T) (Y^*)^2 + c_1|\alpha|_T^2- 4\int_0^{ \tau} v^\prime(Y_s)Z_s dW_s.
\]
This estimate implies that for any $p\geq 1$,
\begin{align}
\mathbb{E}\Big[\Big(\int_0^{\tau}  |Z_s|^2 ds \Big)^{\frac{p}{2}}\Big] \leq c_2 \mathbb{E}\Big[(Y^*)^p + |\alpha|_T^p + \Big|\int_0^{ \tau} v^\prime (Y_s) Z_s dW_s\Big|^{\frac{p}{2}}\Big], \label{lpestimate:1}
\end{align}
where $c_2 := 3^{\frac{p}{2}}\big((3c_1+c^2_1 T) \vee 4\big)^{\frac{p}{2}}$.
Define for each $n\in\mathbb{N}^+$, 
$
\tau_n := \inf \big\{  
t\geq 0:\int_0^t |Z_s|^2 ds
\geq n
\big\} \wedge T.
$
 We then replace $\tau$ by $\tau_n$ and use  Davis-Burkholder-Gundy inequality to obtain
\begin{align*}
c_2\mathbb{E}\Big[ \Big(\int_0^{ \tau_n}v^\prime(Y_s)Z_s dW_s \Big)^\frac{p}{2}\Big]& \leq c_2c({p})M^p
\mathbb{E}\Big[ \Big(\int_0^{ \tau_n} |Y_s|^2 |Z_s|^2 ds\Big)^\frac{p}{4}\Big] \\
& \leq \frac{1}{2}c^2_2c(p)^2 M^{2p}\mathbb{E} \big[ (Y^*)^p\big]  + \frac{1}{2} \mathbb{E}\Big [\Big(\int_0^{\tau_n}  |Z_s|^2 ds\Big)^{\frac{p}{2}}\Big]
\\
&< +\infty.
\end{align*}
We explain that in this inequality, $c(p)$ comes from  Davis-Burkholder-Gundy inequality and only depends on $p$. With this estimate, 
we  come back to (\ref{lpestimate:1}).
Transferring the quadratic term  to the left-hand side  of (\ref{lpestimate:1}) and using Fatou's lemma, we obtain
\[
 \mathbb{E}\Big [\Big(\int_0^{T}  |Z_s|^2 ds\Big)^{\frac{p}{2}} \Big] \leq c \Big( \mathbb{E} \big[ (Y^*)^p+|\alpha|_T^{p}\big] \Big),
\]
where $c : = c_2^2c(p)^2M^{2p} + 2c_2.$

To estimate  $\int_0^T f(|Y_s|)|Z_s|^2 ds$ we use $u: = u^{2f(|\cdot|)}$. This helps to transfer $\int_0^T f(|Y_s|)|Z_s|^2 ds$ to the left-hand side so that standard estimates can be used. 
The proof is omitted since it is not relevant to our study. 
\qed
\end{Proof}

We continue our study 
by sharpening Lemma \ref{lpestimate1} for $p>1$. 
We follow  Proposition 3.2, Briand et al \cite{B2003}  and extend it to quadratic BSDEs. As an important byproduct, we obtain 
the a priori bound for solutions which is crucial to the construction of a solution. 
\begin{Lemma}[A Priori Estimate (ii)]
\label{lpestimate2}
 Let $p > 1$ and   {\rm \ref{lpas1}}  hold for $(F, \xi)$. If $(Y, Z) \in \mathcal{S}^p \times \mathcal{M} $ is a solution to $(F, \xi)$, then	 
\[
\mathbb{E}\big[(Y^*)^p\big] + \mathbb{E}\Big[ \Big( \int_0^T |Z_s|^2 ds   \Big)^{\frac{p}{2}}\Big]
+ \mathbb{E}\Big[ \Big( \int_0^T f(|Y_s|)|Z_s|^2 ds   \Big)^{p}\Big]
  \leq c \Big(\mathbb{E}\big[|\xi|^p + |\alpha|_T^p\big] \Big).
\]
In particular,  
\[
\mathbb{E}\Big[ \sup_{s \in [t, T]}|Y_s|^p \Big|\mathcal{F}_t \Big] \leq c \mathbb{E} \big[|\xi|^p + |\alpha|_{t, T}^p \big| \mathcal{F}_t\big].
\]
In both cases,  $c$ is a constant  only depending on
$T, M^{f(|\cdot|)}, \beta, \gamma, p$.
\end{Lemma} 
\begin{Proof}
Let $u: =u^{f(|\cdot|)}$ and $M : = M^{f(|\cdot|)}$, and  denote 
 $u(|Y_t|), u^\prime(|Y_t|), u^{\prime\prime}(|Y_t|)$ by $u_t, u_t^\prime, u_t^{\prime\prime}$, respectively.
By Tanaka's formula applied  to $|Y_t|$ and  It\^{o}-Krylov formula  applied to $u_t$,
\begin{align*}
 u_t &= u_T + \int_t^{T}\sgn(Y_s) u^\prime_sF(s, Y_s, Z_s)ds -\frac{1}{2}\int_t^{T} \mathbb{I}_{\{Y_s\neq 0\}} u^{\prime\prime}_s|Z_s|^2 ds \\
&- \int_t^{ T} \sgn (Y_s) u^\prime_s Z_s dW_s - \int_t^{ T}u^{\prime}_s dL_s^0(Y),
\end{align*}
where $L^0(Y)$ is the local time of $Y$ at $0$.  Lemma \ref{lpinequality} applied to $u_t$ then gives 
\begin{align*}
|u_t|^p & + \frac{p(p-1)}{2} \int_t^T 
 \mathbb{I}_{\{ u_s \neq 0  \}} 
\mathbb{I}_{\{Y_s \neq 0 \}}|u_s|^{p-2}|u^\prime_s|^2 |Z_s|^2  ds \\
& = |u_T|^p + p \int_t^T \sgn(u_s)|u_s|^{p-1}  \Big (  \sgn(Y_s)u^\prime_sF(s, Y_s, Z_s) -\frac{1}{2} \mathbb{I}_{\{Y_s\neq 0\}}u^{\prime\prime}_s |Z_s|^2\Big)ds\\ &-p\int_t^T\sgn(u_s)|u_s|^{p-1} u^{\prime}_sdL^0_s(Y) 
 - p\int_t^T  \sgn(u_s)\sgn (Y_s)|u_s|^{p-1}u^\prime_s Z_s dW_s.
\end{align*}
To simplify this equality, we recall that $\sgn(u_s)= \mathbb{I}_{\{ u_s \neq 0\}} =\mathbb{I}_{\{ Y_s\neq 0\}}$  and 
$u^{\prime\prime}(x) = 2f(x)u^\prime(x)
$ a.e. Hence
\begin{align*}
|u_t|^p  &+ \frac{p(p-1)}{2} \int_t^T \mathbb{I}_{\{Y_s \neq 0 \}}|u_s|^{p-2}|u^\prime_s|^2 |Z_s|^2  ds \\
& \leq |u_T|^p + p \int_t^T \mathbb{I}_{\{Y_s \neq 0 \}}|u_s|^{p-1} u^\prime_s \big ( \alpha_s + \beta |Y_s| + \gamma |Z_s|\big) ds \\
& - p\int_t^T \sgn (Y_s)|u_s|^{p-1} u^\prime_s Z_s dW_s.
\end{align*}
Let $\{c_n \}_{n\in\mathbb{N}^+}$ be constants to be determined. 
 Since $\frac{|x|}{M}\leq u(|x|)\leq M|x|$ and $\frac{1}{M}\leq  u^\prime(|x|)\leq M$, this inequality yields
\begin{align}
|Y_t|^p  &+ c_1 \int_t^T \mathbb{I}_{\{Y_s \neq 0\}}|Y_s|^{p-2} |Z_s|^2 ds  \nonumber \\
&  \leq M^p |\xi|^p + M^p \int_t^T \mathbb{I}_{\{Y_s \neq 0\}}|Y_s|^{p-1} \big ( \alpha_s + |\beta| |Y_s| + \gamma |Z_s|\big)  ds \nonumber \\
&- p \int_t^T \sgn (Y_s)|u_s|^{p-1}  u^\prime_s Z_s dW_s,\label{up1}
\end{align}
where $c_1 := \frac{p(p-1)}{2M^p} >0$. Observe that in (\ref{up1}),
\[
M^p  \gamma \mathbb{I}_{\{Y_s \neq 0\}}|Y_s|^{p-1}|Z_s| \leq  \frac{M^{2p}\gamma^2}{2c_1}|Y_s|^p + \frac{c_1}{2}\mathbb{I}_{\{Y_s\neq 0\}}|Y_s|^{p-2} |Z_s|^2.
\]
We then use this inequality  to (\ref{up1}).  Set 
 $c_2: = M^p \vee \big( M^p |\beta| + \frac{M^{2p}\gamma^2}{2c_1}\big)$, 
\[
X:  = c_2 \Big( |\xi|^p + \int_0^T |Y_s|^{p-1} \big ( \alpha_s +  |Y_s| \big)  ds\Big),
\]
and $N$ to be the local martingale part of (\ref{up1}). Hence (\ref{up1}) gives
\begin{align}
|Y_t|^p +\frac{c_1}{2} \int_t^T \mathbb{I}_{\{Y_s \neq 0\}}|Y_s|^{p-2} |Z_s|^2 ds \leq X - N_T+ N_t. \label{lp2}
\end{align}
We claim that $N$ is a martingale.
Let $c(1)$ be the constant in Davis-Burkholder-Gundy inequality for $p=1$.  We have
\begin{align*}
\mathbb{E} \big[N^*\big] &\leq c(1) \mathbb{E} \big[\langle N \rangle _T^{\frac{1}{2}}\big] \\
&\leq c(1)M^p \mathbb{E} \Big[ \Big(\int_0^T  |Y_s|^{2p-2} |Z_s|^2 ds \Big)^{\frac{1}{2}}\Big]  \\
& \leq  \frac{c(1)M^p}{p} \Big((p-1)\mathbb{E}\big[(Y^*)^p\big] +  \mathbb{E}\Big[ \Big( \int_0^T |Z_s|^2 ds \Big)^{\frac{p}{2}}\Big]\Big) \\
&< +\infty,
\end{align*}
where the last two lines come from  Young's inequality and Lemma \ref{lpestimate1} (a priori estimate (i)).
 Hence $N$ is a martingale. 
 Coming back to (\ref{lp2}), we deduce that 
\begin{align}
\mathbb{E}\Big[  \int_0^T \mathbb{I}_{\{Y_s\neq 0\}}|Y_s|^{p-2} |Z_s|^2 ds \Big] \leq \frac{2}{c_1} \mathbb{E} [X]. \label{2x}
\end{align}
Now we estimate $Y$ via $X$. To this end,  taking supremum over $t\in [0, T]$ and using Davis-Burkholder-Gundy inequality to (\ref{lp2})  give
\begin{align}
\mathbb{E}\big[(Y^*)^p\big] \leq \mathbb{E}[ X] + c(1) \mathbb{E} \big[\langle N \rangle_T^\frac{1}{2}\big]. \label{lp3}
\end{align}
The second term  in (\ref{lp3}) yields by Cauchy-Schwartz inequality that 
\begin{align*}
 c(1)\mathbb{E} [\langle N \rangle_T^\frac{1}{2}] & \leq  c(1)M^p \mathbb{E}\Big[ (Y^*)^{\frac{p}{2}} \Big( \int_0^T\mathbb{I}_{\{Y_s \neq 0\}} |Y_s|^{p-2} |Z_s|^2  ds\Big)^{\frac{1}{2}}\Big] \\
 &\leq \frac{1}{2}\mathbb{E} \big[(Y^*)^p\big] + \frac{c(1)^2M^{2p}}{2} \mathbb{E}\Big[\int_0^T \mathbb{I}_{\{Y_s \neq 0\}}|Y_s|^{p-2}|Z_s|^2  ds\Big].
\end{align*}
Using   (\ref{2x}) to this inequality gives the estimate of
$\langle N\rangle^{\frac{1}{2}}$  via $Y$ and $X$.   With this estimate we come back to 
 (\ref{lp3}) and obtain
\[
\mathbb{E}[(Y^*)^p] \leq 2\Big(1+\frac{2c({1})^2M^{2p}}{c_1}\Big) \mathbb{E}[X].
\]
Set $c_3: = 2c_2\big( 1+ \frac{c(1)^2M^{2p}}{2}\big)  .$
This inequality  yields
\begin{align}
\mathbb{E}[(Y^*)^p] \leq c_{3} \Big(\mathbb{E}\big[|\xi|^p \big] + \mathbb{E}\Big[ \int_0^T |Y_s|^{p-1}\alpha_s ds\Big]+ \mathbb{E}\Big[ \int_0^T |Y_s|^p ds \Big]\Big).\label{up2}
\end{align} 
 Young's inequality used to the second term on the right-hand side of this inequality gives
\[
c_{3} \int_0^T |Y_s|^{p-1}  \alpha_s ds \leq \frac{1}{2} (Y^*)^p + \frac{c_3}{p}\Big(\frac{2}{c_3 q}\Big)^{\frac{p}{q}}|\alpha|_T^p,
\]
where $q$ is the conjugate index of $p$.
Set $c_4: = 2\Big(c_3 \vee\frac{c_3}{p}\big(\frac{2}{c_3 q}\big)^{\frac{p}{q}}  \Big)$.  (\ref{up2}) and the above inequality  yield
\[
\mathbb{E}\big[(Y^*)^p\big] \leq  c_4\Big( \mathbb{E}\big[|\xi|^p + |\alpha|_T^p\big]+  \mathbb{E}\Big[ \int_0^T \sup_{u\in[0, s]} |Y_u|^p ds \Big]\Big),
\]
By Gronwall's lemma, 
\[
\mathbb{E}\big[(Y^*)^p\big] \leq c_4 \exp(c_4 T) \mathbb{E}\big[|\xi|^p + |\alpha|_T^p\big].
\]
Finally, by Lemma \ref{lpestimate1} we conclude that there exists a constant $c$ only depending on $T, M, \beta, \gamma, p$ such that 
\[
\mathbb{E}\big[(Y^*)^p\big] + \mathbb{E}\Big[ \Big( \int_0^T |Z_s|^2 ds   \Big)^{\frac{p}{2}}\Big]
+ \mathbb{E}\Big[ \Big( \int_0^T f(|Y_s|)|Z_s|^2 ds   \Big)^{p}\Big]
  \leq c \mathbb{E}\big[|\xi|^p + |\alpha|_T^p\big] .
\]
To prove the remaining statement, we view any fixed $t\in [0, T]$ as  the initial time, reset
\[
X:  = c_2\Big( |\xi|^p +  \int_t^T |Y_s|^{p-1} \big ( \alpha_s +  |Y_s| \big)  ds\Big)
\]
and replace all estimates by conditional estimates. 
\qed
\end{Proof}

An immediate consequence of Lemma \ref{lpestimate2} is that 
\[
|Y_t| \leq \Big(c \mathbb{E} \big[|\xi|^p + |\alpha|_{T}^p \big| \mathcal{F}_t\big] \Big)^{\frac{1}{p}},
\]
i.e., $Y$  has an a priori bound  which is a  continuous supermartingale. 

With this estimate    
we are  ready to construct a $\mathbb{L}^p(p>1)$ solution
via inf-(sup-)convolution
as in  Briand et al   \cite{BH2006}, \cite{BH2008},  \cite{BLS2007}.
A localization procedure where the a priori bound plays a crucial role is used  and the monotone stability  result takes the limit. 
\begin{thm}[Existence]\label{lpexistence} Let $p > 1$ and {\rm\ref{lpas1}} hold for $(F, \xi)$.
Then there exists a solution of $(F, \xi)$ in $\mathcal{S}^p\times \mathcal{M}^p$.
\end{thm}
\begin{Proof}
We introduce the notations used throughout the proof. Define the process
\[
X_t := \Big(c \mathbb{E}\big[|\xi|^p + |\alpha|_T^p \big|\mathcal{F}_t \big]\Big)^{\frac{1}{p}},
\]
where $c$ is the constant defined in Lemma \ref{lpestimate2}. Obviously
$X$ is continuous by It\^{o}  representation theorem. Moreover, for  each $m, n \in \mathbb{N}^+$, set 
\begin{align*}
\tau_m &: = \inf \big\{
t\geq 0: |\alpha|_t +X_t \geq m      \big\} \wedge T, \\
\sigma_n &:= \inf \big\{ 
t\geq 0: |\alpha|_t \geq n\big\} \wedge T.
\end{align*}
 It then follows from the continuity of $|\alpha|_\cdot$ and $X$ that $\tau_m$ and $\sigma_n$ increase stationarily to $T$ as $m, n$ goes to $+\infty$, respectively. To apply a double approximation procedure, we  define
\begin{align*}
F^{n, k}(t, y, z): &=\mathbb{I}_{\{t\leq \sigma_n \}} \inf_{y^\prime, z^\prime}\big \{ F^+(t, y^\prime, z^\prime) + n|y-y^\prime|+n|z-z^\prime| \big\} \\
&-\mathbb{I}_{\{t\leq \sigma_k \}} \inf_{y^\prime, z^\prime} \big\{ F^-(t, y^\prime, z^\prime) + k|y-y^\prime|+k|z-z^\prime| \big\},
\end{align*}
and $
\xi^{n , k}: = \xi^+ \wedge n - \xi^- \wedge k.$ 

Before proceeding to the proof we give some useful facts.  
By  Lepeltier and San Martin \cite{LS1997}, 
$F^{n, k}$ is  Lipschitz-continuous in $(y, z)$; as  $k$ goes to $+\infty$,  $F^{n, k}$ converges decreasingly 
 uniformly on compact sets to a limit denoted by $F^{n,\infty}$; as $n$ goes to $+\infty$, $F^{n, \infty}$ converges increasingly uniformly on compact sets to $F$. Moreover, $\big||F^{n, k}(\cdot, 0, 0) |\big|_T$ and $\xi^{n, k}$ are bounded. 
 
 Hence,  by  Briand et al \cite{B2003},
there exists a unique solution $(Y^{n, k}, Z^{n, k}) \in \mathcal{S}^p\times\mathcal{M}^{p}$ of $(F^{n, k},  \xi^{n, k})$; by comparison theorem, $Y^{n, k}$ is increasing in $n$ and decreasing in $k$. We are about to take the limit by the  monotone stability result.

However,  $\big||F^{n, k}(\cdot, 0, 0)|\big|_T$ and $Y^{n, k}$
are not uniformly bounded in general. 
To overcome this difficulty, we  use Lemma \ref{lpestimate2}
and  work on random time interval
where $Y^{n,k}$ and $\big||F^{n, k}(\cdot, 0, 0)|\big|_\cdot$  are uniformly bounded.
This is the motivation to 
 introduce $X$ and $\tau_m$. To be more precise, the localization procedure is as follows. 

Note that $(F^{n, k}, \xi^{n, k})$ satisfies  \ref{lpas1} associated with $(\alpha, \beta, \gamma, \varphi, f)$. Hence by  Lemma \ref{lpestimate2} (a priori estimate (ii)),
\begin{align}
|Y_t^{n, k}| &\leq \Big(c \mathbb{E}\big[ |\xi^{n, k}|^p + |\mathbb{I}_{[0, \sigma_n \vee \sigma_k]}\alpha|_T^p \big| \mathcal{F}_t\big]\Big)^{\frac{1}{p}}\nonumber\\
&\leq 
X_t. \label{lpx1}
\end{align}
In view of the definition of $\tau_m$,  we deduce that 
\begin{align}
|Y^{n, k}_{t\wedge \tau_m} |\leq X_{t\wedge \tau_m} \leq m.\label{lpbound}
\end{align}
Hence $Y^{n, k}$ is uniformly bounded on $[0, \tau_m]$. Secondly, 
 given $(Y^{n, k}, Z^{n, k})$ which solves $(F^{n, k}, \xi^{n, k})$, it is  immediate  that 
 $(Y^{n, k}_{\cdot\wedge\tau_m}, \mathbb{I}_{[0, \tau_m]}Z^{n, k})$ solves 
$(\mathbb{I}_{[0, \tau_m]}F^{n, k}, Y^{n, k}_{\tau_m})$. 
To make the monotone stability result adaptable, we 
use a truncation procedure. Define 
\[
\rho(y): = -\mathbb{I}_{\{y<-m\}} m + \mathbb{I}_{\{|y| \leq m\}} y +\mathbb{I}_{\{y>m\}}  m.
\]
Hence from (\ref{lpbound})  $(Y^{n, k}_{\cdot\wedge\tau_m}, \mathbb{I}_{[0, \tau_m]}Z^{n, k})$ meanwhile solves $(\mathbb{I}_{[0, \tau_m]}(t)F^{n, k}(t, \rho(y), z), Y^{n, k}_{{\tau_m}})$. Secondly,  we have
\begin{align*}
|\mathbb{I}_{[0, \tau_m]}(t)F^{ n, k} (t, \rho(y), z)|&\leq \mathbb{I}_{\{t\leq \tau_m\}} \big(\alpha_t +  \varphi(|\rho(y)|) + \gamma |z| + f(|\rho(y)|)|z|^2\big)\\
&\leq\mathbb{I}_{\{t\leq \tau_m\}} \Big( \alpha_t+ \varphi(m) + \gamma |z|  + \sup_{|y|\leq m}f(|\rho(y)|) |z|^2\Big)\\
&\leq \mathbb{I}_{\{t\leq \tau_m\}} \Big(\alpha_t+\varphi (m) + \frac{\gamma^2}{4} +  \big(\sup_{|y|\leq m}f(|\rho(y)|) +1\big) |z|^2\Big),
\end{align*}
where $\sup_{|y|\leq m}f(|\rho(y)|)$ is bounded  for each $m$ due to $f(|\cdot|)\in\mathcal{I}$. Moreover, the definition of $\tau_m$ implies $|\alpha|_{\tau_m} \leq m$.
Hence we can use the monotone stability result (see Kobylanksi \cite{K2000} or Briand and  Hu \cite{BH2008}) 
to obtain   $(Y^{m, n, \infty}, Z^{ m, n, \infty})\in\mathcal{S}^\infty\times\mathcal{M}^2$ which solves $(\mathbb{I}_{[0, \tau_m]}(t)F^{n, \infty}(t, \rho(y), z), \inf_k Y_{\tau_m}^{n, k})$. Moreover,   $Y^{m, n, \infty}_{\cdot \wedge \tau_m}$ is the $\mathbb{P}$-a.s. uniform limit of $Y^{n, k}_{\cdot \wedge \tau_m}$ as $k$ goes to $+\infty$.
These arguments hold for any $m, n\in\mathbb{N}^+$.  

Due to this convergence result we can pass the comparison property to 
 $Y^{m, n, \infty}$. We use the monotone stability result again to the sequence indexed by $n$ to obtain $(\widetilde{Y}^m, \widetilde{Z}^m) \in \mathcal{S}^\infty\times\mathcal{M}^2$ which solves $(\mathbb{I}_{[0, \tau_m]}(t)F(t, \rho(y), z), \ 
\sup_n\inf_k Y^{n, k}_{\tau_m})$.   Likewise, $\widetilde{Y}^m_{\cdot}$ is the $\mathbb{P}$-a.s. uniform limit of $Y^{m, n, \infty}_{\cdot}$  as $n$ goes to $+\infty$. Hence we  obtain from (\ref{lpbound}) that $|\widetilde{Y}^m_{t}| \leq X_{t\wedge \tau_m} \leq m $. Therefore, 
$(\widetilde{Y}^m, \widetilde{Z}^m)$  solves $(\mathbb{I}_{[0, \tau_m]}F, \sup_n\inf_k Y_{\tau_m}^{n, k})$, i.e., 
\begin{align}
\widetilde{Y}^m_{t\wedge \tau_m}= \sup_n\inf_k {Y}_{\tau_m}^{n, k} +
\int_{t\wedge \tau_m}^{\tau_m} F(s, \widetilde{Y}^m_s, \widetilde{Z}_s^m)ds -
\int_{t\wedge \tau_m}^{\tau_m}  \widetilde{Z}_s^m dW_s. \label{locbsde}
\end{align}

We recall that the monotone stability result also implies that 
 $\widetilde{Z}^m$ is the $\mathcal{M}^2$-limit of $\mathbb{I}_{[0, \tau_m]}Z^{n, k}$ as $k, n$ goes to $+\infty$. This fact and previous convergence results give
 \begin{align}
 \widetilde{Y}^{m+1}_{\cdot\wedge \tau_m}  &=
  \widetilde{Y}^{m}_{\cdot\wedge \tau_m} \ \mathbb{P}\text{-a.s.}, \nonumber\\
 \mathbb{I}_{\{t\leq \tau_m\}}\widetilde{Z}^{m+1}_t &= \mathbb{I}_{\{t\leq \tau_m\}}\widetilde{Z}^m_t \ dt\otimes d\mathbb{P}\text{-a.e.} \label{piece}
 \end{align} 
Define $(Y, Z)$ on $[0, T]$ by
\begin{align*}
Y_t := \mathbb{I}_{\{t\leq \tau_1\}} \widetilde{Y}_t^1  + \sum_{m\geq 2} \mathbb{I}_{]\tau_{m-1}, \tau_m]}\widetilde{Y}^m_t, \\
Z_t := \mathbb{I}_{\{t\leq \tau_1\}} \widetilde{Z}_t^1  + \sum_{m\geq 2} \mathbb{I}_{]\tau_{m-1}, \tau_m]}\widetilde{Z}^m_t.
\end{align*}
By (\ref{piece}),  we have $Y_{\cdot \wedge \tau_m} = \widetilde{Y}^m_{\cdot\wedge \tau_m}$
and
$\mathbb{I}_{\{t\leq \tau_m\}}Z_t = \mathbb{I}_{\{t\leq \tau_m\}} \widetilde{Z}_t^m$. Hence we can rewrite  
 (\ref{locbsde}) as
\[
Y_{t\wedge \tau_m} = \sup_n\inf_k Y^{n, k}_{\tau_m} +\int_{t\wedge \tau_m}^{\tau_m}  F{(s, Y_s, Z_s )}ds -\int_{t\wedge \tau_m}^{\tau_m} Z_s dW_s.
\]
 By sending $m$ to $+\infty$, we deduce  that $(Y, Z)$ solves $(F, \xi)$.  Since $(Y^{n, k}, Z^{n, k})$ verifies Lemma \ref{lpestimate2}, we can use Fatou's lemma to prove that $(Y, Z)\in\mathcal{S}^p\times\mathcal{M}^p$.
\qed
\end{Proof}

Theorem \ref{lpexistence} proves the existence of a $\mathbb{L}^p(p>1)$ solution under  \ref{lpas1} which to our knowledge the most general asssumption. 
For example, 
 \ref{lpas1}(ii) allows one to get rid of monotonicity in $y$ which is required by, e.g.,  Pardoux \cite{P1999} and Briand et al \cite{BC2000}, \cite{B2003}, \cite{BLS2007}. Meanwhile, in contrast to these works, 
the generator can also be quadratic by setting $f(|\cdot|)\in\mathcal{I}$. 
Hence Theorem \ref{lpexistence} provides a unified way to construct solutions of both non-quadratic and quadratic BSDEs via the monotone stability result.

 On the other hand, Theorem \ref{lpexistence} is an extension of  Bahlali et al  \cite{BEO2014} which only studies BSDEs with $\mathbb{L}^2$ integrability and linear-quadratic growth. However, in contrast to their work,  \ref{lpas1} is not sufficient in our setting to ensure the existence of a maximal or minimal solution, since the  double approximation procedure  makes the comparison between solutions impossible.
 
However, to prove the existence of a maximal or minimal solution is no way impossible. Since we have $X$ as the a priori bound for solutions, 
we can convert the question of existence into the question of existence  for quadratic BSDEs with double  barriers. This problem has been solved by introducing the notion of generalized BSDEs; see Essaky and Hassani \cite{E2013}.
 
Let us turn to the uniqueness result. Motivated by Briand and  Hu \cite{BH2008} or  Da Lio and  Ley \cite{LL2006} from the point of view of PDEs, we impose a convexity condition so as to use $\theta$-techinique which proves to be convenient to treat quadratic generators.  
We start from  comparison theorem and then move to uniqueness and stability result.
To this end, the following assumptions on $(F, \xi)$ are needed. 
\begin{as}\label{lpas3}  Let $p > 1.$
There exist $\beta_1, \beta_2 \in \mathbb{R}$, $\gamma_1, \gamma_2 \geq 0$, an $\mathbb{R}^+$-valued $\Prog$-measurable process $\alpha$,  a continuous nondecreasing function $\varphi: \mathbb{R}^+ \rightarrow \mathbb{R}^+$ with $\varphi(0)=0$, 
$f(|\cdot|)\in\mathcal{I}$  and $F_1, F_2: \Omega\times[0, T]\times \mathbb{R}\times\mathbb{R}^d\rightarrow \mathbb{R}$ 
which are $\Prog\otimes\mathcal{B}(\mathbb{R})\otimes\mathcal{B}(\mathbb{R}^d)$-measurable
such that $F= F_1+ F_2$, 
$|\xi| + |\alpha|_T\in\mathbb{L}^p$
and  $\mathbb{P}$-a.s.
 \begin{enumerate}
 \item [(i)] for any $t\in[0, T]$, 
$(y, z)\longmapsto F(t, y, z)$ is continuous; 
\item[(ii)] $F_1(t, y, z)$ is  monotonic in $y$ and Lipschitz-continuous in $z$,  and $F_2(t, y, z)$ is monotonic at $y=0$ and of linear-quadratic growth in $z$, 
 i.e., for any  $t\in [0, T], y, y^\prime \in \mathbb{R}, z, z^\prime \in \mathbb{R}^d$,
\begin{align*}
\sgn (y -y ^\prime)\big(F_1(t, y, z)- F_1(t, y^\prime, z)\big)&\leq \beta_1 |y- y^\prime|, \\
 \big|F_1(t, y, z)- F_1(t, y, z^\prime)\big| &\leq  \gamma_1 |z-z^\prime|, \\
 \sgn(y)F_2(t, y, z) &\leq \beta_2|y|+ \gamma_2 |z| + f(|y|)|z|^2;
\end{align*}
\item [(iii)] for any $t\in[0, T]$, $(y, z)\longmapsto F_2(t, y, z)$  is convex;
\item [(iv)] for any $(t, y, z)\in [0, T]\times\mathbb{R}\times\mathbb{R}^d$, 
\begin{align*}
&|F(t, y, z)|\leq \alpha_t + \varphi(|y|) + (\gamma_1 +\gamma_2)|z|+ f(|y|)|z|^2.
  \end{align*}
\end{enumerate}
\end{as}

Intuitively, \ref{lpas3} specifies an additive structure consisting of two  classes of BSDEs. 
The cases where $F_2 = 0$  coincide with 
classic existence and uniqueness results for $\mathbb{R}$-valued BSDEs; see, e.g.,  Pardoux \cite{PP1990} or Briand et al \cite{BC2000}, \cite{B2003}. 
When $F_1 =0$,  the BSDEs concern and generalize those studied by
 Bahlali et al \cite{BEK2013}. 
Given  convexity as an additional requirement, we can prove an existence and uniqueness result in the presence of both components. This can be seen as a general version of  the additive structure discussed in Section \ref{degenerate} and 
a complement to  the quadratic BSDEs studied 
by Bahlali et al \cite{BEK2013} and
Briand and  Hu \cite{BH2008}.

We start our proof of comparison theorem by observing that  \ref{lpas3} implies \ref{lpas1}. Hence the  existence of a
$\mathbb{L}^p(p>1)$  solution is ensured. 

\begin{thm}[Comparison]\label{lpcompare} Let $p > 1$,
and  $(Y, Z), (Y^\prime, Z^\prime) \in \mathcal{S}^p\times\mathcal{M}$ be solutions of $(F, \xi)$, $(F^\prime, \xi^\prime)$, respectively. If  $\mathbb{P}$-a.s. for any $(t ,y ,z)\in [0, T]\times\mathbb{R}\times\mathbb{R}^d$, 
$F(t, y, z)\leq F^\prime(t, y, z)$ and $\xi \leq \xi^\prime$, and $F$ verifies {\rm\ref{lpas3}},
then $\mathbb{P}$-a.s.
$Y_\cdot\leq Y^\prime_\cdot$.
\end{thm}
\begin{Proof}
We introduce the  notations used throughout the proof. 
For any $\theta \in (0, 1)$,
define
\begin{align*}
\delta F_t &:= F(t, Y_t^\prime, Z_t^\prime) - F^\prime(t, Y_t^\prime, Z_t^\prime), \\
\delta_\theta Y&: = Y- \theta Y^\prime,  \\
\delta Y &: = Y- Y^\prime, 
\end{align*} 
and  $\delta_\theta Z, \delta Z$, etc. analogously.
$\theta$-technique applied to the generators yields
\begin{align}
& F(t, Y_t, Z_t) - \theta F^\prime(t, Y_t^\prime, Z_t^\prime) \nonumber\\
&=\big(F(t, Y_t, Z_t) - \theta F(t, Y^\prime_t, Z^\prime_t)\big)+\theta\big(F(t, Y^\prime_t, Z^\prime_t) -F^\prime(t, Y_t^\prime, Z_t^\prime)\big)  \nonumber\\
&= \theta \delta F_t + \big(F(t, Y_t, Z_t) -\theta F(t,  Y_t^\prime, Z_t^\prime)\big) \nonumber \\
& = \theta \delta F_t + \big(F_1(t, Y_t, Z_t\big) -\theta F_1(t, Y_t^\prime, Z_t^\prime)\big) + \big(F_2(t, Y_t, Z_t) -\theta F_2(t, Y_t^\prime, Z_t^\prime)\big).
 \label{lpconvex1}
\end{align}
By  \ref{lpas3}(iii),
\begin{align*}
F_2(t, Y_t, Z_t) &= F_2(t, \theta Y_t^\prime + (1-\theta)\frac{\delta_\theta Y_t}{1-\theta}, \theta Z_t^\prime + (1-\theta)\frac{\delta_\theta Z_t}{1-\theta})\nonumber\\
& \leq \theta F_2(t, Y_t^\prime, Z_t^\prime) + (1-\theta) F_2(t, \frac{\delta_\theta Y_t}{1-\theta}, \frac{\delta_\theta Z_t}{1-\theta}). 
\end{align*}
Hence we have 
\begin{align}
F_2(t, Y_t, Z_t) -  \theta F_2(t, Y_t^\prime, Z_t^\prime) \leq (1-\theta) F_2(t, \frac{\delta_\theta Y_t}{1-\theta}, \frac{\delta_\theta Z_t}{1-\theta}).
\label{lpconvex2}
\end{align}
Let $u$  be the function defined in Section \ref{I} associated with a function of class $\mathcal{I}$ to be determined later. 
 Denote  $u((\delta_\theta Y_t)^+), u^\prime((\delta_\theta Y_t )^+), u^{\prime\prime}((\delta_\theta Y_t)^+)$ by
 $u_t, u_t^\prime$, $u_t^{\prime\prime}$, respectively. It is then known from Section \ref{I} that  $u_t \geq 0$ and $u^\prime_t >0$. 
 For any $\tau \in \mathcal{T}$, 
 Tanaka's formula applied  to $(\delta_\theta Y)^+$, It\^{o}-Krylov formula applied  to $u((\delta_\theta Y_t)^+)$ and Lemma \ref{lpinequality} give
\begin{align}
|u_{t\wedge \tau}|^p  &+ \frac{p(p-1)}{2} \int_{t\wedge \tau}^\tau \mathbb{I}_{\{ \delta_\theta Y_s>0\}}|u_s|^{p-2}|u_s^\prime|^2  |\delta_\theta Z_s|^2 ds \nonumber\\
& \leq |u_\tau |^p + p \int_{t\wedge \tau}^\tau \mathbb{I}_{\{ \delta_\theta Y_s > 0\}}|u_s|^{p-1} \underbrace{\Big(u^\prime_s\big ( F(s, Y_s, Z_s)- \theta F^\prime (s, Y^\prime_s, Z^\prime_s)\big) -\frac{1}{2}u^{\prime\prime}_s |\delta_\theta Z_s|^2 \Big)}_{: =\Delta_s}ds \nonumber \\
& 
- p\int_{t\wedge \tau}^\tau \mathbb{I}_{\{\delta_\theta Y_s >0\}}|u_s|^{p-1} u^\prime_s  \delta_\theta Z_s  dW_s. \label{lpconvex3}
\end{align}
By (\ref{lpconvex1}), (\ref{lpconvex2}),  \ref{lpas3}(ii) and   $\delta F \leq 0$,  we deduce that, on  $\{ \delta_\theta Y_s >0\}$, 
\begin{align*}
 \Delta_s
&\leq u^\prime_s \Big(  F_1 (s, Y_s, Z_s) - \theta F_1 (s, Y_s^\prime, Z_s^\prime) +\beta_2 (\delta_\theta Y_s)^+ + \gamma_2 |\delta_\theta Z_s| +\frac{f(\frac{|\delta_\theta Y_s|}{1-\theta})}{1-\theta}|\delta_\theta Z_s|^2\Big)
-\frac{1}{2}u_s^{\prime\prime} |\delta_\theta Z_s|^2.
\end{align*}  To eliminate the quadratic term, 
we associate $u$ with 
$\frac{f(\frac{|\cdot|}{1-\theta})}{1-\theta}$, i.e.,  
\begin{align*}
u(x): &= \int_0^x \exp\Big (2 \int_0^y \frac{f(\frac{|u|}{1-\theta})}{1-\theta}du\Big)dy \\
&= \int_0^x \exp\Big (2 \int_0^{\frac{y}{1-\theta}}f(|u|)du\Big)dy. 
\end{align*}
Hence,  on $\{ \delta_\theta Y_s >0\}$, the above inequality  gives
\begin{align}
 \Delta_s
 \leq  u^\prime_s \big(
  F_1 (s, Y_s, Z_s) - \theta F_1 (s, Y_s^\prime, Z_s^\prime) +\beta_2 (\delta_\theta Y_s)^+ + \gamma_2 |\delta_\theta Z_s|\big).
  \label{lpconvex4}
\end{align}
We are about to send $\theta$ to $1$, and to this end we give some auxiliary facts. Reset
$
M: = \exp \big(2\int_0^\infty f(u)du\big).
$
Obviously
 $1\leq M< +\infty$.
By dominated convergence, for $x\geq 0,$ we have
\begin{align}
&\lim_{\theta\rightarrow 1} u(x)= M x, \nonumber\\
&\lim_{\theta\rightarrow 1} u^\prime(x) = M \mathbb{I}_{\{ x>0\}} + \mathbb{I}_{\{ x= 0\}}. \label{lim}
\end{align}
Taking (\ref{lpconvex4}) and (\ref{lim}) into account, we  come back to  (\ref{lpconvex3}) and 
 send $\theta$ to $1$. 
  Fatou's lemma used to the $ds$-integral on the left-hand side of  (\ref{lpconvex3}) and 
dominated convergence used to the rest integrals give
\begin{align}
&((\delta Y_{t\wedge \tau})^+)^p + 
\frac{p(p-1)}{2}  \int_{t\wedge \tau}^\tau \mathbb{I}_{\{\delta Y_s >0\}} ((\delta Y_s)^+)^{p-2}|\delta Z_s|^2ds \nonumber\\
  &\leq   ((\delta Y_\tau)^+)^p + p \int_{t\wedge \tau}^\tau \mathbb{I}_{\{ \delta Y_s >0\}}((\delta Y_s)^+)^{p-1} ( 
F_1(s, Y_s, Z_s) -F_1(s, Y_s^\prime, Z_s^\prime) +   
 \beta_2 (\delta Y_s)^+ + \gamma_2 |\delta Z_s|     )  ds \nonumber \\
&
-p\int_{t\wedge \tau}^\tau \mathbb{I}_{\{\delta Y_s >0\}} ((\delta Y_s)^+)^{p-1}\delta Z_s  dW_s. \label{lpconvex5}
\end{align}
 Moreover,  
\ref{lpas3}(ii) implies 
\[
\mathbb{I}_{\{ \delta Y_s > 0\}}\big(F_1(s, Y_s, Z_s) - F_1(s, Y_s^\prime, Z_s^\prime) \big)\leq 
\mathbb{I}_{\{ \delta Y_s > 0\}}\big(\beta_1 (\delta Y_s)^+ + \gamma_1 |\delta Z_s|\big).
\]
We then use this inequality to (\ref{lpconvex5}). To eliminate the local martingale,
we replace $\tau$  by its localization sequence $\{\tau_n\}_{n\in\mathbb{N}^+}$.
By the same way of estimation as  in  Lemma \ref{lpestimate2} (a priori estimate (ii)), we obtain
\[
 ((\delta Y_{t\wedge \tau_n})^+ )^p 
  \leq c\mathbb{E}  \big[ ((\delta Y_{\tau_n})^+ )^p \big| \mathcal{F}_t\big] ,
\]
where $c$ is a constant  only depending on $T, \beta_1, \beta_2, \gamma_1, \gamma_2, p$. Since $Y, Y^\prime \in \mathcal{S}^p$ and  $\mathbb{P}$-a.s. $\xi \leq \xi^\prime$, dominated convergence  yields  $\mathbb{P}$-a.s. $Y_t\leq Y_t^\prime$.
Finally, by the continuity of $Y$ and  $Y^\prime$ we conclude that  $\mathbb{P}$-a.s. $Y_\cdot \leq Y^\prime_\cdot$.
\qed
\end{Proof}

As a byproduct, we obtain the following existence and uniqueness result.
\begin{Corollary}[Uniqueness]\label{lpunique} Let {\rm\ref{lpas3}} hold for $(F, \xi)$. Then there exists a unique solution in $\mathcal{S}^p\times \mathcal{M}^p$.
\end{Corollary}
\begin{Proof}
\ref{lpas3} implies \ref{lpas1}. Hence existence result holds. The uniqueness is immediate from   Theorem \ref{lpcompare} (comparison theorem).
\qed
\end{Proof}

It turns out that a 
 stability result also holds given the convexity condition. We  denote $(F, \xi)$ satisfying \ref{lpas3} by
$(F, F_1, F_2, \xi).$ We set $\mathbb{N}^0:= \mathbb{N}^+ \cup \{ 0 \}.$
\begin{prop}[Stability]
\label{lpstability}Let $p> 1$. Let
$(F^n, F_1^n, F_2^n, \xi^n)_{n\in\mathbb{N}^0}$ satisfy {\rm\ref{lpas3}}  associated with 
$(\alpha^n, \beta_1, \beta_2, \gamma_1, \gamma_2, \varphi, f)$, and $(Y^n, Z^n)$ be their unique solutions in $\mathcal{S}^p\times\mathcal{M}^p$, respectively. 
If 
$
\xi^n - \xi {\longrightarrow} 0$ and $\int_0^T |F^n -F^0 |(s, Y_s^0, Z_s^0) ds {\longrightarrow} 0$
in $\mathbb{L}^p$
as $n$ goes to $+\infty$, 
then $(Y^n, Z^n)$ converges to $(Y, Z)$ in $\mathcal{S}^p\times\mathcal{M}^p$.
\end{prop}
\begin{Proof}
We prove the stability result in the spirit of  Theorem \ref{lpcompare} (comparison theorem).
 For any $\theta \in (0, 1)$,   define 
\begin{align*}
\delta F_t^n &:= F^0(t, Y^0_t, Z^0_t) - F^n(t, Y^0_t, Z^0_t), \\
\delta_\theta Y^n&: = Y^0-  \theta Y^n, \\
\delta Y^n&: = Y^0-   Y^n,
\end{align*} 
and  $\delta_\theta Z^n, \delta Z^n$, etc. analogously.
We  observe the $\theta$-difference of the generators. Likewise, \ref{lpas3}(iii) implies that 
\begin{align*}
&F^0(t, Y^0_t, Z^0_t) - \theta F^n (t, Y^n_t, Z^n_t)\\
&= \delta F_t^n + \big(  F^{n}(t, Y^0_t, Z^0_t) -\theta F^n(t, Y^n_t, Z^n_t)    \big)\\
&\leq 
\delta F_t^n + \big(F_1^n (t, Y_t^0, Z_t^0) - \theta F_1^n(t, Y_t^n, Z_t^n)\big) +
(1-\theta) F^n_2 (t, \frac{\delta_\theta Y_s^n}{1-\theta}, \frac{\delta_\theta Z_s^n}{1-\theta}).
\end{align*}
We first prove convergence of $Y^n$ and later use it to show  that $Z^n$ also converges.

(i).  
By exactly the same arguments as in Theorem 
\ref{lpcompare} but keeping $\delta F_t^n$ along the deductions,
we obtain
\begin{align}
&((\delta Y_t^n)^+)^p + 
\frac{p(p-1)}{2} \int_t^T\mathbb{I}_{\{ \delta Y_s^n >0\}}  ((\delta Y_s^n)^+)^{p-2}|\delta Z_s^n|^2ds \nonumber\\
  &\leq  ((\delta\xi^n)^+)^p+  p \int_t^T\mathbb{I}_{\{\delta Y_s^n >0\}}  ((\delta Y_s^n)^+)^{p-1} \big(  |\delta F_s^n| + (\beta_1 + \beta_2)(\delta Y_s^n)^+ + (\gamma_1+\gamma_2)|\delta Z_s^n|     \big) ds \nonumber\\
&
-p\int_t^T \mathbb{I}_{\{\delta Y_s^n >0\}}((\delta Y_s^n)^+)^{p-1}\delta Z_s^n  dW_s. \label{convex4}
\end{align}
 By the same way of estimation as in  Lemma \ref{lpestimate2} (a priori estimate (ii)),  we obtain
\[
\mathbb{E}\big[\big(((\delta Y^n)^+)^*\big)^p\big] \leq c\Big( \mathbb{E} \big[\big((\delta \xi^n)^+\big)^p\big] +  \mathbb{E}\big[ \big||\delta F^n_\cdot|\big|_T^p \big] \Big),
\] 
where $c$ is a constant only depending on
$T, \beta_1, \beta_2, \gamma_1, \gamma_2, p$.
Interchanging $Y^0$ and $Y^n$ and analogous deductions then yield
\[
\mathbb{E}\big[\big(((-\delta Y^n)^+)^*\big)^p\big]  \leq c\Big( \mathbb{E} \big[\big((-\delta \xi^n)^+\big)^p\big] +  \mathbb{E}\big[ \big||\delta F^n_\cdot|\big|_T^p \big] \Big).
\] 
Hence a combination of the two inequalities implies the convergence of $Y^n.$

(ii).
To prove the convergence of $Z^n$, we combine the arguments in
Lemma \ref{lpestimate1} (a priori estimate (i)) and Theorem \ref{lpcompare}. 
To this end, we introduce the function $v$ defined in Section \ref{I} associated with 
a function of class $\mathcal{I}$ to be determined later. By It\^o-Krylov formula,	
\begin{align}
v(\delta_\theta Y_0^n) &= v(\delta_\theta \xi^n)
+\int_0^T  v^\prime(\delta_\theta Y_s^n) \big(F^0(s, Y_s^0, Z_s^0) -\theta F^n(s, Y_s^n, Z_s^n) \big) ds
\nonumber\\ 
&-\frac{1}{2}\int_0^T v^{\prime\prime}(\delta_\theta Y_s^n)|\delta_\theta Z_s^n|^2 ds 
- \int_0^T v^\prime(\delta_\theta Y_s^n)\delta_\theta Z_s^n dW_s. \label{lpstability1}
\end{align} Note that \ref{lpas3}(ii)(iii) and 
$v^\prime(\delta_\theta Y_s^n) = \sgn(\delta_\theta Y_s^n)|v^\prime(\delta_\theta Y_s^n)|$ give 
\begin{align}
v^\prime(\delta_\theta Y_s^n) \big(F^0(s, Y_s^0, Z_s^0) &-\theta F^n(s, Y_s^n, Z_s^n) \big) \nonumber\\
&\leq
|v^\prime(\delta_\theta Y_s^n)||\delta F_s^n| \nonumber\\
&+ |v^\prime(\delta_\theta Y_s^n)|\sgn(\delta_\theta Y_s^n)\big(
F^n_1(s, Y_s^0, Z_s^0) - \theta F^n_1(s, Y_s^n, Z_s^n)
\big)
\nonumber\\
&+ |v^\prime(\delta_\theta Y_s^n)|
\Big(  \beta_2 |\delta_\theta Y_s^n| 
+\gamma_2 |\delta_\theta Z_s^n|
+ \frac{f(\frac{|\delta_\theta Y_s^n|}{1-\theta})}{1-\theta} |\delta_\theta Z_s^n|^2
\Big). \label{vf}
\end{align}
We associate $v$ with $\frac{f(\frac{|\cdot|}{1-\theta})}{1-\theta}$ so as to eliminate the quadratic term. Note that
\begin{align}
&\lim_{\theta \rightarrow 1} v(x) =\frac{1}{2}|x|^2, 
\nonumber\\
&\lim_{\theta \rightarrow 1} v^\prime(x) =x. \label{limv}
\end{align}
With (\ref{vf}), (\ref{limv}) and \ref{lpas3}(ii), we come back to (\ref{lpstability1})  and send $\theta$ to $1$.   This gives
\begin{align*}
\frac{1}{2}\int_0^T |\delta Z_s^n|^2 ds
&\leq \frac{1}{2}|\delta \xi^n|^2
+\int_0^T |\delta Y_s^n| \big(|\delta F_s^n| + 
(|\beta_1|+|\beta_2|) |\delta Y_s^n| + (\gamma_1 +\gamma_2)|\delta Z_s^n|
 \big) ds\\
 &- \int_0^T  \delta Y_s^n\delta Z_s^n dW_s.
\end{align*}
Now we use the same way of estimation as in Lemma \ref{lpestimate1} to obtain
\begin{align*}
\mathbb{E}\Big[\Big( \int_0^T  |\delta Z_s^n|^2ds\Big)^{\frac{p}{2}}\Big]
\leq 
c \mathbb{E}\big[
((\delta Y^n)^*)^p +
\big||\delta F^n_\cdot|\big|_T^p
\big
],
\end{align*}
where $c$ is a constant only depending on $T, \beta_1, \beta_2, \gamma_1, \gamma_2, p$.
The convergence of $Z^n$ is then immediate from (i).
\qed
\end{Proof}
\begin{remm}
So far we have obtained the  existence and uniqueness of  a  $\mathbb{L}^p(p > 1)$ solution. The solvability for $p =1$ is not included due to the failure of  Lemma \ref{lpestimate2} (a priori estimate (ii)). One may overcome this difficulty by imposing additional structure conditions as in  Briand et al \cite{B2003}, \cite{BH2006}. To save pages the analysis of $\mathbb{L}^1$ solutions is hence omitted. 
\end{remm}

\subsection{Applications to Quadratic PDEs}
\label{qpde}
In this section, we give an application of our results  to quadratic PDEs. More precisely, we prove the probablistic representation for  the nonlinear  Feymann-Kac formula associated with the BSDEs in our study. Let us consider the following semilinear PDE
\begin{align}
&\partial_t u(t, x) +\mathcal{L}u(t, x) + 
F(t, x, u(t, x), \sigma^\top \nabla_x u(t, x)) =0, \nonumber\\
&u(T, \cdot) = g, \label{pde}
\end{align}
where $\mathcal{L}$ is the infinitesimal generator of the solution  $X^{t_0, x_0}$ to the Markovian SDE
\begin{align}
X_t &= x_0 +\int_{t_0}^t b(s, X_s)ds + \int_{t_0}^t \sigma (s, X_s)dB_s, \label{sde}
\end{align}
for any $(t_0, x_0)\in [0, T]\times \mathbb{R}^n$, 
$t\in [t_0, T]$.  
Denote 
  a solution to the BSDE 
\begin{align}
Y_t = g(X_T^{t_0, x_0})
+\int_t^T F(s, X_s^{t_0, x_0}, Y_s, Z_s)ds -\int_t^T Z_sdW_s, \ t\in[t_0, T],\label{fbsde}
\end{align}
by  $(Y^{t_0, x_0}, Z^{t_0, x_0})$ or $(Y, Z)$ when there is no ambiguity. 
The probablistic representation for nonlinear Feymann-Kac formula consists of proving that,  in Markovian setting, $u(t, x):=Y_t^{t, x}$
is a solution at least in the viscosity sense to  (\ref{pde}) when the source of nonlinearity $F$ is quadratic in $\nabla_x u(t, x)$ and  $g$ is an unbounded function.  To put it more precisely, let us  introduce the FBSDEs.

{\bf The Forward Markovian SDEs.}
Let $b: [0, T]\times\mathbb{R}^n \rightarrow \mathbb{R}^n$, $\sigma: [0, T]\times\mathbb{R}^d\rightarrow \mathbb{R}^{n\times d}$ be continuous functions and assume there exists $\beta \geq 0$ such that $\mathbb{P}$-a.s. for any $t\in[0, T]$, $|b(t, 0)| +|\sigma(t, 0)| \leq \beta$ and $b(t, x), \sigma(t, x)$ are  Lipschitz-continuous in $x$, i.e., $\mathbb{P}$-a.s.  for any $t\in[0, T]$, $x, x^\prime \in\mathbb{R}^n$, 
\[
|b(t, x)- b(t, x^\prime)|
+|\sigma(t, x)- \sigma(t, x^\prime)|\leq \beta |x-x^\prime|.
\]
Then for any $(t_0, x_0)\in [0, T]\times\mathbb{R}^n$,  (\ref{sde}) has a unique solution $X^{t_0, x_0}$ in  $\mathcal{S}^p$ for any $p \geq 1$.  

{\bf The Markovian BSDE.} We continue with the setting of the forward equations above. Set $q\geq 1$.
Let $F_1, F_2: [0, T]\times\mathbb{R}^n\times\mathbb{R}\times\mathbb{R}^d \rightarrow \mathbb{R}$, 
$g: \mathbb{R}^n\rightarrow \mathbb{R}$ be continuous functions, $\varphi: \mathbb{R}^+ \rightarrow \mathbb{R}^+$  a continuous nondecreasing function with $\varphi(0)=0$ and
$f(|\cdot|)\in \mathcal{I}$, and assume moreover  $F = F_1 + F_2$ such  that
\begin{enumerate}
\item[(i)] 
 $F_1(t, x, y, z)$ is  monotonic in $y$ and  Lipschitz-continuous in $z$, and $F_2(t, x, y, z)$ is monotonic at $y=0$ and of linear-quadratic growth in $z$, 
 i.e., for any $(t, x)\in[0, T]\times\mathbb{R}^n$, $y, y^\prime \in \mathbb{R}, z, z^\prime \in \mathbb{R}^d$,
\begin{align*}
\sgn (y -y ^\prime)\big(F_1(t, x, y, z)- F_1(t,  x, y^\prime, z)\big)&\leq \beta |y- y^\prime|, \\
 \big|F_1(t, x, y, z)- F_1(t, x, y, z^\prime)\big| &\leq  \beta |z-z^\prime|, \\
 \sgn(y)F_2(t, x, y, z) &\leq \beta|y| + \beta |z| + f(|y|)|z|^2;
\end{align*}
\item[(ii)] $(y, z) \longmapsto F_2(t, x, y, z)$ is  convex ;
\item[(iii)] for any $(t, x, y, z)\in [0, T]\times\mathbb{R}^n\times\mathbb{R}\times\mathbb{R}^d$, 
\begin{align*}
| F (t, x, y, z) | &\leq  \beta \big(1+ |x|^q + 2|z|\big) + \varphi(|y|)
  + f(|y|)|z|^2, \\
 |g(x)| &\leq \beta  \big(1 + |x|^q\big). 
  \end{align*}
\end{enumerate} 

Since $X^{t_0, x_0}\in \mathcal{S}^p$ for any $p \geq 1$, the  above structure conditions on  $F$ and $g$ allow one to use Corollary \ref{lpunique}
to construct a unique solution $(Y^{t_0, x_0}, Z^{t_0, x_0})$ in $\mathcal{S}^p\times\mathcal{M}^p$ of 
 (\ref{fbsde})
 for any $p > 1$.
Moreover, by standard arguments,
$Y_{t_0}^{t_0, x_0}$
 is deterministic for any $(t_0, x_0) \in [0, T]\times\mathbb{R}^n$. Hence $u(t, x)$ defined as $Y_t^{t, x}$ is a deterministic function. With this fact 
 we now turn to the main result of this section: $u$ is a viscosity solution of   (\ref{pde}). Before our proof let us recall the definition of a viscosity solution. 
 
{\bf Viscosity Solution.} A continuous function $u: [0, T]\times\mathbb{R}^n \rightarrow \mathbb{R}$ is called a viscosity subsolution (respectively supersolution) to (\ref{pde}) if   $u(T, x) \leq g(x)$ (respectively $u(T, x) \geq g(x)$)
and 
for any smooth function  $\phi$ such that  $u-\phi$ reaches the local maximum (respectively local minimum) at $(t_0, x_0)$, we have 
\[
\partial_t \phi(t_0, x_0)
+\mathcal{L}\phi(t_0, x_0)
+ F(t_0, x_0, u(t_0, x_0), \sigma^\top \nabla_x \phi(t_0, x_0))\geq 0\ (\text{respectively} \leq 0).
\]
A function $u$ is called a viscosity solution to (\ref{pde}) if it is both a viscosity subsolution and supersolution. 
\begin{prop}
Given the above assumptions, $u(t, x)$  is continuous with 
\[
|u(t, x)| \leq c\big( 1+ |x|^q\big),
\] where $c$ is a constant. Moreover, $u$ is a viscosity solution to {\rm(\ref{pde})}.
\end{prop}
\begin{Proof}
 Due to the Lipschitz-continuity of $b$ and $\sigma$, $X^{t, x}$ is continuous in $(t, x)$, e.g.,  in mean square sense. The continuity of $u$ is then an immediate consequence of  Theorem \ref{lpstability} (stability). The proof relies on standard arguments and  hence is omitted.
 By  Lemma \ref{lpestimate2} (a priori estimate (ii)), we prove that $u$ satisfies the above polynomial growth. It thus remains to prove that $u$ is a viscosity solution to (\ref{pde}).

Let $\phi$ be a smooth function such that $u-\phi$ reaches local maximum at $(t_0, x_0)$. Without loss of generality we  assume that the local maximum is global and 
$u(t_0, x_0) = \phi(t_0, x_0).$ We aim at proving 
\[
\partial_t \phi(t_0, x_0)
+\mathcal{L}\phi(t_0, x_0)
+ F(t_0, x_0, u(t_0, x_0), \sigma^\top \nabla_x \phi(t_0, x_0))\geq 0.
\]
From (\ref{fbsde}) we obtain
\[
Y_t = Y_{t_0}
-\int_{t_0}^t F(s, X_s^{t_0, x_0}, Y_s, Z_s)ds +\int_{t_0}^t Z_sdW_s. \]
By It\^{o}'s formula, 
\[
\phi(t, X_t^{t_0, x_0}) = 
\phi(t_0, x_0)
+
\int_{t_0}^{t}
\big\{\partial_s\phi+\mathcal{L}\phi \big\} (s, X_s^{t_0, x_0})ds +
\int_{t_0}^{t}
\sigma^{\top} \nabla_x \phi(s, X_s^{t_0, x_0}) dW_s.  
\]
Now we take any $t\in [t_0, T]$.
Note that the existence of a unique solution of (\ref{sde}) and (\ref{fbsde})  implies by Markov property 
that $Y_t =u(t, X_t^{t_0, x_0})$.  Hence,
 $\phi(t, X_t^{t_0, x_0})\geq u(t, X_t^{t_0, x_0}) =    Y_t$.  By touching property, on the set
$
\big\{ \phi(t, X_t^{t_0, x_0}) = Y_t\big\}
$ we have
\begin{align*}
&\partial_t \phi(t, X_t^{t_0, x_0}) +\mathcal{L}\phi(t, X_t^{t_0, x_0}) + F(t, X_t^{t_0, x_0}, Y_t, Z_t) \geq 0\ \ \mathbb{P}\text{-a.s.}, \\
&\sigma^\top \nabla_x\phi(t, X_t^{t_0, x_0}) -Z_t =0\ \  \mathbb{P}\text{-a.s.}
\end{align*}
Now we set $t=t_0$. We have
$\phi(t_0, X_{t_0}^{t_0, x_0}) = \phi(t_0, x_0) = u(t_0, x_0) = Y_{t_0}$. Moreover, 
the above equality implies 
$Z_{t_0} = \sigma^\top \nabla_x \phi(t_0, x_0)$.  Plugging the two equalities into the above inequality gives
\[
\partial_t \phi(t_0, x_0)
+\mathcal{L}\phi(t_0, x_0)
+ F(t_0, x_0, u(t_0, x_0), \sigma^\top \nabla_x \phi(t_0, x_0))\geq 0.
\]
 Hence   $u$ is a viscosity subsolution of (\ref{pde}).  $u$ being a viscosity supersolution and thus a viscosity solution can be proved analogously.
\qed
\end{Proof}
\textbf{Acknowledgement.} The author thanks Martin Schweizer for his supervision and many helpful remarks.
\bibliography{myrefs}
\bibliographystyle{plain}
\end{document}